\documentclass[12pt,reqno]{amsart}

\usepackage{amssymb}
\usepackage{amsthm}
\usepackage{latexsym,amsmath,amssymb}
\usepackage{ascmac} 
\usepackage[pdftex]{graphicx}
\usepackage{braket}
\usepackage{amscd}
\usepackage{amsfonts}
\usepackage{mathrsfs}
\usepackage{bm}  
\usepackage{enumitem}  
\usepackage{hhline} 
\usepackage{appendix}
\usepackage[colorlinks=true, bookmarks=true,
bookmarksnumbered=true, bookmarkstype=toc, linkcolor=blue,
urlcolor=blue, citecolor=blue]{hyperref}

\usepackage{float}
\usepackage{wrapfig}
\restylefloat{figure}
\usepackage{floatflt}

\def\XXint#1#2#3{{\setbox0=\hbox{$#1{#2#3}{\int}$}
\vcenter{\hbox{$#2#3$}}\kern-.5\wd0}}

\theoremstyle{plain} 
\newtheorem{thm}{Theorem}[section]
\newtheorem{lem}[thm]{Lemma}
\newtheorem{cor}[thm]{Corollary}
\newtheorem{prop}[thm]{Proposition}

\newtheorem{example}[thm]{Example}
\theoremstyle{definition} 
\newtheorem{dfn}[thm]{Definition}
\newtheorem{rem}[thm]{Remark}
\makeatletter
    
    \@addtoreset{equation}{section}
  \makeatother

\makeatletter
\def\tr{\mathop{\operator@font tr}\nolimits}  
\makeatother

\makeatletter
\def\dist{\mathop{\operator@font dist}\nolimits}  
\makeatother

\makeatletter
\def\div{\mathop{\operator@font div}\nolimits}  
\makeatother

\makeatletter
\def\exp{\mathop{\operator@font exp}\nolimits}  
\makeatother

\makeatletter
\def\essinf{\mathop{\operator@font {\itshape ess}.\inf}\nolimits}  
\makeatother

\makeatletter
\def\esssup{\mathop{\operator@font {\itshape ess.}\sup}\nolimits}  
\makeatother

\newcommand{\R}{\mathbb{R}}
\newcommand{\N}{\mathbb{N}}

\newcommand{\M}{{\mathcal M}}
\newcommand{\A}{{\mathcal A}}

\def\phi{\varphi}

\def\ol{\overline}

\begin{document}

\title[Harnck inequalities and H\"older estimates with weak scaling conditions]{Harnck inequalities and H\"older estimates for fully nonlinear integro-differential equations with weak scaling conditions}

\author[S. KITANO]{
Shuhei Kitano
}


\address{
Department of Applied Physics\endgraf
Waseda University\endgraf
Tokyo, 169-8555\endgraf
JAPAN
}
\email{sk.koryo@moegi.waseda.jp}

\maketitle

\begin{abstract}
H\"older estimates and Harnack inequalities are studied for fully nonlinear integro-differential equations under some mild assumptions. We allow the kernels  of variable order and critically close to 2.
\end{abstract}

\section{Introduction}
\label{sec:intro}

In this paper, we study the fully nonlinear nonlocal equations of the form:
\begin{align}
Iu(x):=\sup_{a\in\mathcal{A}}\inf_{b\in\mathcal{B}}\int_{\R^n}\delta(u,x,y)\kappa_{a,b}(x,y)\nu(x,|y|)
dy=f(x)\quad\mbox{in }B_R,\label{eq}
\end{align}
where $u:\R^n\to\R$ is an unknown function, $f:\R^n\to\R$ is a given function and $B_R$ is an open ball with its center $0$ and its radius $R>0$.
We write $\delta(u,x,y):=u(x+y)+u(x-y)-2u(x)$.
Let $0<\lambda\leq\Lambda$. Suppose that for any $a\in\mathcal{A}$, $b\in\mathcal{B}$ and $x\in B_R$, $y\in\R^n$, $\kappa_{a,b}(x,y)$ satisfies $\kappa_{a,b}(x,y)=\kappa_{a,b}(x,-y)$ and
\begin{equation}\label{UE}
\lambda\leq \kappa_{a,b}(x,y) \leq \Lambda.
\end{equation}
$Iu$ is well defined for $C^{1,1}$ and bounded functions $u$ if $\nu:B_R\times (0,\infty)\to(0,\infty)$ satisfies
\begin{equation}\label{eq12}
\int_{\R^n}(1\wedge |y|^2)\nu(x,|y|)dy<\infty\quad \text{for }x\in B_R,
\end{equation}
where we denote $l\wedge m:= \min\{l,m\}$ for $l,m\in\R^n$.

The main purpose of this paper is to present the H\"older estimates and the Harnack inequalities for solutions of \eqref{eq} under the weak scaling condition of $\nu$, which we introduce below. 
When $\nu(x,r)=r^{-n-\sigma}$ for $\sigma\in(0,2)$, the first result about Harnack inequalties and H\"older estimates was established in \cite{BL02}, by using a probabilistic approach
and the H\"older estimate was obtained in \cite{Sil06}, where the proof is based on the theory of partial differential equations.
Here we note that those results are nonlocal versions of a classical result of~\cite{KS} by Krylov and Safonov for non-divergence second order equations, but blow up as the order $\sigma$ of the equation approaches to $2$.
Harnack inequalities and H\"older estimates whose constants do not blow up as $\sigma \to 2$ were obtained in \cite{CS09}.

There are more general kernels, for example
\begin{equation}\label{eq11}
\nu(x,r)\asymp \frac{1}{r^{n+\sigma}}(\log(2/r))^\gamma
\end{equation}
for $0<r<1$, $\sigma\in[0,2]$ and $\gamma\in\R$.
Here, we used the notation $f(r)\asymp g(r)$ which means that the quotient $f(r)/g(r)$ stays between two positive constants.
Note that when $\sigma=2$, we need $\gamma<-1$ so that \eqref{eq12} holds.
If $\sigma\in(0,2)$, the H\"older estimate is obtained by~\cite{BK05b, KM13, KM17, KKL16, KL, Sil06} and the Harnack inequality is also proved by~\cite{BK05a, KM13, KKL16, KL}.
An interesting case is when $\sigma=0$ and in this case, some scale invariant results fail in general, but the paper~\cite{KM17} proved the equi-continuity of solutions, which is not covered in this paper.
On the other hand, all of these previous results cannot be applied when $\sigma=2$ and $\gamma<-1$, although these kernels are treated in~\cite{Grz14,Mim13}, in the case where the nonlocal operator $I$ is linear and independent of $x$ i.e. $I$ is of form
\[
Iu(x)=\int_{\R^n}\delta(u,x,y)\nu(|y|)dy.
\]
One of our main contribution is that we provide the Harnack inequality and the H\"older estimate for fully nonlinear nonlocal equations with the kernels of form \eqref{eq11} when $\sigma=2$ and $\gamma<-1$.

We are also interested in variable order kernel which is of form:
\[
\nu(x,r)\asymp \frac{1}{r^{n+\sigma(x)}}
\]
for $\sigma(x)\in(0,2)$.
In this case, the Harnack inequality were given by \cite{BK05a} and the H\"older estimate, by~\cite{BK05b, Sil06}.
Here the Harnack inequality of~\cite{BK05a} is not scale invariant (in other words, its constant depends on the radial $R$ of the domain).
In the case of nonlocal Dirichlet form or divergence form equations with variable order kernels, the scale invariant Harnack inequality was established by~\cite{CKW}.
Another main contribution of this paper is that we show the scale invariant Harnack inequality for fully nonlinear integral equations with variable order kernel under some scale invariant conditions.

It can be possible to extent our results to the non-symmetric operator $I$ with an anisotropic kernel $\nu$, in the sense that $I$ is of form
\[
Iu(x)=\sup_{a\in\mathcal{A}}\inf_{b\in\mathcal{B}}\int_{\R^n}(u(x+y)-u(x)-\nabla u(x)\cdot y\chi_{B_1}(y))\kappa_{a,b}(x,y)\nu(x,y)dy,
\]
where we do not assume the symmetry $\kappa(x,y)=\kappa(x,-y)$ and the ellipticity condition \eqref{UE} can be relaxed.
In this direction, \cite{CLD} studied non-symmetric operators and \cite{BCF, KRS}, more general conditions of \eqref{UE} when $\nu(x,y)\asymp |y|^{-n-\sigma}$.
Also, anisotropic kernels were treated in~\cite{BK05a, BK05b, KM13}.
It is worth to mention that for those operators, the Harnack inequality fails in general even when the H\"older estimate holds and several counterexamples were given by~\cite{BBCK09, BK05a, BS05, BS07, F09a}.
Finally we note that there are many related work, which we did not mention above;
see, for instance \cite{B88, BKK, C20, CCW, CZ16, DK20, F09a, F09, GS17, GS18, KL2} and references therein.




\subsection{Assumptions and Main results}
\begin{dfn}
Set
\begin{align*}
h(x,r):=\int_{\R^n}\left(1\wedge\frac{|y|^2}{r^2}\right)\nu(x,|y|)dy.
\end{align*}
We say $\nu$ satisfies (A) if there exist $c_U,\alpha>0$ such that for $x\in B_R$,
\begin{align}
h(x,1)&=1,\tag{A1}\label{A1}\\
0\leq \nu(x,r)&\leq\nu(x,s)\quad\mbox{for }0<s\leq r\mbox{ and}\tag{A2}\label{A2}\\
h(x,tr)&\leq c_Ut^{-\alpha}h(x,r)\quad\mbox{for }r>0\mbox{ and }t\geq1.\tag{A3}\label{A3}
\end{align}
\end{dfn}
\noindent
Especially, we call \eqref{A3} the {\it weak scaling condition} by following \cite{GS17,GS18}.
In Section \ref{sec2}, we will recall some basic properties of $h$ and equivalent conditions with the weak scaling condition, which were investigated by \cite{GS17,GS18}.

One of our main result is the following H\"older estimate:
\begin{thm}\label{main thm1}There are positive constants $\eta\in(0.1)$ and $C\geq1$ depending only on $n,\ \lambda,\ \Lambda,\ c_U$ and $\alpha$  with the following property. Let $\nu$ be a function satisfying (A). If $u\in C(B_R)\cap L^\infty(\R^n)$ is a viscosity  solution of $Iu=f$ in $B_R$, then it follows that
\begin{equation*}
\| u\|_{C^\eta(B_{R/2})}\leq CR^{-\eta}(\|u\|_{L^\infty(\R^n)}+\|h(\cdot,R)^{-1}f(\cdot)\|_{L^\infty(B_R)}).
\end{equation*}
\end{thm}

We also obtain the Harnack inequality. In this case, we assume in addition the following condition: there exist $c^*\geq1$ such that
\begin{align}
\frac{\nu(x_1,r)}{h(x_1,s)}\leq c^*\frac{\nu(x_2,r)}{h(x_2,s)}\quad\mbox{for }|x_1-x_2|\leq s \leq r\label{B}\tag{B}
\end{align}

\begin{thm}\label{mainthm2}
There exists a constant $C=C(n,\lambda,\Lambda,c_U,\alpha,c^*)>0$ with the following property. Let $\nu$ satisfy (A) and \eqref{B}. If $u\in C(B_{2\sqrt{n}R})\cap L^\infty(\R^n)$ is non-negative in $\R^n$ and  a viscosity solution of $Iu=f$ in $B_{2\sqrt{n}R}$, then it follows that
\[
\sup_{B_{R/2}}u\leq C\left(\inf_{\mathcal{Q}_{3R}}u+ \|h(\cdot,R)^{-1}f\|_{L^\infty(B_{2\sqrt{n}R})}\right).
\]
\end{thm}
\begin{rem}
We fill a gap in the assumptions of \cite{Grz14}. 
In the case of the nonlocal operator $I$ is linear and constant coefficient and $n\geq 3$, the condition (A) is equivalent with the one in \cite{Grz14} (and \eqref{B} is not needed in this case since $I$ is independent of $x$).
On the other hand, when $n=1,2$, additional conditions are required in \cite{Grz14} although it is not necessary in our case (see Corollary 6 and Theorem 7 in \cite{Grz14}).
\end{rem}
\begin{rem}
The condition~\eqref{B} is a modification of assumptions in \cite{BK05a, F09a, CKW}, and cannot be avoided for the Harnack inequality to hold.
In Section \ref{sect8}, we present a counterexample to the Harnack inequality when we assume only (A) but not \eqref{B}.
There is a novelty because our example is isotropic, whereas all examples obtained in~\cite{BBCK09, BK05a, BS05, BS07, F09a} are anisotropic.
\end{rem}

\subsection{Examples}
Here we provide some examples of kernels $\nu_i$ for $i=1,...,5$.
We denote $h_i(x,r):=\int_{\R^n}\left(1\wedge\frac{|y|^2}{r^2}\right)\nu_i(x,|y|)dy$ in what follow.
\begin{example}[Fractional laplacian]\label{ex16}
For any $\sigma\in [\alpha,2)$, let 
\[
\nu_1(x,r):=\A(n,\sigma) r^{-n-\sigma}\quad\mbox{for }r>0,
\]
where $\A(n,\sigma)$ is the normalized constant defined by
\begin{align*}
\A(n,\sigma)
&:=\left(\int_{\R^n}(1\wedge|y|^2)|y|^{-n-\sigma}dy
\right)^{-1}\\
&=\left(\int_0^1\int_{\partial B_1}r^{1-\sigma}d\omega dr+\int_{1}^\infty\int_{\partial B_1} r^{-1-\sigma}d\omega dr
\right)^{-1}\\
&=\left(\frac{|\partial B_1|}{(2-\sigma)}+\frac{|\partial B_1|}{\sigma}\right)^{-1}
=\frac{(2-\sigma)\sigma}{2|\partial B_1|}.
\end{align*}
Then, we have $h_1(x,r)=r^{-\sigma}$ and $\nu_1$ satisfies the condition (A).
In fact, it is easy to see \eqref{A1} and \eqref{A2} hold and we also have
\[
h_1(x,tr)=(tr)^{-\sigma}\leq t^{-\alpha}r^{-\sigma}=t^{-\alpha}h(x,r)
\]
for $t\geq1$, $r>0$ and $\alpha\in(0,\sigma]$.
\end{example}
\begin{example}[Sum of fractional laplacians]
For $\alpha>0$, let $\sigma_i\in[\alpha,2)$ for $i=1,2,...$ and $\{a_i\}_{i=1}^\infty$ be a sequence of positive constants such that $\sum_{i=1}^\infty a_i=1$.
Then,
\[
\nu_2(x,r):=\sum_{i=1}^\infty a_i\A(n,\sigma_i)r^{-n-\sigma_i}\quad\mbox{for }r>0.
\]
and $h_2(x,r)=\sum_{i=1}^\infty a_i r^{-\sigma_i}$ satisfy the condition (A), which is obtained in a similar way with Example 1.6.
\end{example}
\begin{example}For $\alpha>0$, let $\gamma\in\R$, $\sigma\in[0,2]$ and
\[
\nu_3(x,r):=
\frac{\mathcal{B}(n,\sigma,\gamma)}{r^{n+\sigma}}\left\{
\begin{array}{ll}
\left(1-\log r\right)^{\gamma}&\mbox{for }r\in (0, 1],\\
\left(1+\log r\right)^{-\gamma}&\mbox{for }r>1,
\end{array}
\right.
\]
where we need $m>1$ when $\sigma=0$ and $m<-1$ when $\sigma=2$ so that $(1\wedge|y|^2)\nu_3(x,|y|)$ is integrable and $\mathcal{B}(n,\sigma,\gamma)$ is the normalized constant so that $h_3(x,1)=1$.
For any $\epsilon>0$, taking account of
\begin{align*}
\lim_{r\to+0}\frac{t^{n+\sigma-\epsilon}\nu_3(x,tr)}{\nu(x,r)}
&=\lim_{r\to\infty}\frac{t^{n+\sigma-\epsilon}\nu_3(x,tr)}{\nu(x,r)}=t^{-\epsilon}\leq1
\end{align*}
for $t\geq1$ and
\begin{align*}
\lim_{t\to\infty}\frac{t^{n+\sigma-\epsilon}\nu_3(x,tr)}{\nu(x,r)}
&=0,
\end{align*}
there exists $c_1=c_1(n,\sigma,\gamma,\epsilon)>0$ such that $\nu(x,tr)\leq c_1t^{-n-\sigma+\epsilon}\nu(x,r)$.
Hence, if $\sigma>0$, we can choose $\epsilon<\sigma$ and then we have
\begin{align*}
h_3(x,tr)
&=\int_{\R^n}\left(1\wedge\frac{|y|^2}{(tr)^2}\right)\nu(x,|y|)dy\\
&=\int_{\R^n}\left(1\wedge\frac{|y|^2}{r^2}\right)\nu(x,t|y|)t^ndy\leq C_{n,\sigma,\gamma,\epsilon}t^{-\sigma+\epsilon}h_3(x,r).
\end{align*}
Therefore \eqref{A3} holds with $\alpha=\sigma-\epsilon$ and $c_U=c_1$.
When $\sigma=0$, we can compute
\[
h_3(x,r)\asymp
\left\{
\begin{array}{ll}
\left(1-\log r\right)^{\gamma+1}&\mbox{for }r\in (0, 1],\\
\left(1+\log r\right)^{-\gamma+1}&\mbox{for }r>1.
\end{array}
\right.
\]
Hence \eqref{A3} does not hold for any $\alpha,c_U>0$.
\end{example}   
\begin{example}[Variable order fractional laplacian]\label{ex19}
For $\alpha>0$, consider $\sigma:B_R\to[\alpha,2)$.
Then,
\[
\nu_4(x,r):=\A(n,\sigma(x))r^{-n-\sigma(x)}
\]
satisfies the condition (A).
However the condition \eqref{B} fails unless $\sigma(\cdot)$ is a constant and moreover the Harnack inequality also fails in general.
Let us consider a modified kernel defined by
\[
\nu_5(x,r):=\frac{1}{2}\left(\tilde{\A}(n,\sigma(x))r^{-n-\sigma(x)}\chi_{(0,1)}(r)+\A(n,s_0)r^{-n-s_0}\right)
\]
for $s_0\in[\alpha,2)$ and $\tilde{\A}(n,\sigma(x)):=(\int_{B_1}|y|^{2-n-\sigma(x)}dy)^{-1}=(2-\sigma(x))|\partial B_1|^{-1}$.
Then,
\begin{equation}\label{h5}
h_5(x,r)=\frac{1}{2}\left(\left(\frac{2}{\sigma(x)}r^{-\sigma(x)}-\frac{2-\sigma(x)}{\sigma(x)}\right)\chi_{[0,1)}(r)+r^{-2}\chi_{[1,\infty)}(r)+r^{-s_0}\right).
\end{equation}
\eqref{B} holds if $\nu_5$ satisfies the following condition (i) or (ii):
\[
\text{(i) }\sup_{B_R}\sigma<2 \text{ and }|\sigma(x_1)-\sigma(x_2)|\leq\frac{M_1}{\log(2/|x_1-x_2|)}\quad\mbox{for }|x_1-x_2|<1
\]
for some constant $M_1>0$ or
\[
\text{(ii) }\sup_{B_R}\sigma\leq\gamma_0+s_0 \text{ and }|\sigma(x_1)-\sigma(x_2)|\leq M_2|x_1-x_2|^{\gamma_0}\quad\mbox{for }|x_1-x_2|<1
\]
for some constants $M_2>0$ and $\gamma_0\in(0,1]$.
Indeed, in the case of (i), we have
\begin{align*}
\nu_5(x_1,r)
&=\frac{1}{2}\left(\frac{2-\sigma(x_1)}{2-\sigma(x_2)}\tilde{\A}(n,\sigma(x_2))e^{(\sigma(x_2)-\sigma(x_1))\log r}r^{-n-\sigma(x_2)}\chi_{(0,1)}(r)+\A(n,s_0)r^{-n-s_0}\right)\\
&\leq \frac{1}{2}\left(\frac{2-\alpha}{2-\sup\sigma}e^{M_1\log r/\log(2/r) }\tilde{\A}(n,\sigma(x_2))r^{-n-\sigma(y)}\chi_{(0,1)}(r)+\A(n,s_0)r^{-n-s_0}\right)\\
&\leq c_2 \nu_5(x_2,r)
\end{align*}
for $r\in(0,1)$, $|x_1-x_2|\leq r$ and $c_2:=(2-\alpha)\sup_{0\leq r\leq1}e^{M_1\log r/\log(2/r) }/(2-\sup \sigma)>0$ and we can see $h_5(x_1,r)\leq c_3h_5(x_2,r)$ for $r\in(0,1)$ and some $c_3>0$ by an analogous argument.
Hence we arrive at \eqref{B} for $r\in(0,1)$.
Since \eqref{B} is clear when $r\geq1$, \eqref{B} holds for any $r>0$.
On the other hand, if $\nu_5$ satisfies (ii), we observe
\begin{align*}
&\quad\nu_5(x_1,r)\\
&=\frac{1}{2}\left(\left(\tilde{\A}(n,\sigma(x_2))+\frac{\sigma(x_2)-\sigma(x_1)}{|\partial B_1|}\right)e^{(\sigma(x_2)-\sigma(x_1))\log r}r^{-n-\sigma(x_2)}\chi_{(0,1)}(r)+\A(n,s_0)r^{-n-s_0}\right)\\
&\leq \frac{1}{2}\left(\tilde{\A}(n,\sigma(x_2))r^{-n-\sigma(y)}\chi_{(0,1)}(r)+\frac{1}{|\partial B_1|}r^{-n-\sigma(y)+\gamma_0}+\A(n,s_0)r^{-n-s_0}\right)\\
&\leq \frac{1}{2}\left(\tilde{\A}(n,\sigma(x_2))r^{-n-\sigma(y)}\chi_{(0,1)}(r)+\left(\A(n,s_0)+\frac{1}{|\partial B_1|}\right)r^{-n-s_0}\right)\\
&\leq c_4 \nu_5(x_2,r)
\end{align*}
for $r\in(0,1)$, $|x_1-x_2|\leq r$ and $c_4:=(\A(n,s)+1/|\partial B_1|)/\A(n,s_0)$.
$h_5(x_1,r)\leq c_5 h_5(x_2,r)$ can be also obtained for $r\in(0,1)$ and some $c_5>0$ by a similar way and hence we obtain the condition \eqref{B} in the case of (ii).
\end{example}
\begin{rem}\label{rem110}
Similar conditions for $\sigma(x)$ to (i) of Example \ref{ex19} were considered in~\cite{B88, BK05b, BKK}.
On the other hand, up to our knowledge, the condition (ii) was not treated in the literature.
It allows $\sup_{B_R}\sigma=2$ and in this case, the nonlocal operator may contain second order differential terms as we can observe
\begin{align*}
&\quad I_{\sigma} u(x)\\
&:=\sup_{a\in\mathcal{A}}\inf_{b\in\mathcal{B}}\int_{\R^n}\delta(u,x,y)\kappa_{a,b}(x,y)\frac{1}{2}\left(\frac{\tilde{\A}(n,\sigma)}{|y|^{n+\sigma}}\chi_{(0,1)}(r)+\frac{\A(n,s_0)}{|y|^{n+s_0}}\right)dy\\
&\to \sup_{a\in\mathcal{A}}\inf_{b\in\mathcal{B}}\int_{\partial B_1}\left\langle D^2u(x)\omega,\omega\right\rangle \lim_{h\to0+0}\kappa_{a,b}(x,h\omega)\frac{d\omega}{|\partial B_1|}+\int_{\R^n}\delta(u,x,y)\kappa_{a,b}(x,y)\frac{\A(n,s_0)}{2|y|^{n+s_0}}dy\\
&=:I_{2}u(x),
\end{align*}
as $\sigma\to2$.
Below, we provide a modification of Theorem \ref{mainthm2} for nonlocal equations with second order terms. 
\end{rem}
\begin{cor}
There exists a constant $C=C(n,\lambda,\Lambda,\alpha,\gamma_0,s_0)>0$ with the following property. Let $s_0+\gamma_0\geq2$, $\sigma:B_R\to [\alpha,2]$ be a function satisfying the condition (ii) and $I_{\sigma}u(x)$ be an operator as in Remark \ref{rem110} for $\sigma\in [\alpha,2]$. If $u\in C^2(B_{2\sqrt{n}R})\cap L^\infty(\R^n)$ is non-negative in $\R^n$ and a classical solution of $I_{\sigma(x)}u=f$ in $B_{2\sqrt{n}R}$, then it follows that
\[
\sup_{B_{R/2}}u\leq C\left(\inf_{\mathcal{Q}_{3R}}u+ \|h_5(\cdot,R)^{-1}f\|_{L^\infty(B_{2\sqrt{n}R})}\right),
\]
where $h_5$ is from \eqref{h5}.
\end{cor}
\begin{proof}
For small $\rho>0$, consider the nonlocal operator $I_{(\sigma(x)-\rho)}u(x)$, which still satisfies the condition (A) and (ii) of Example \ref{ex19}.
Note that (ii) implies \eqref{B}.
Hence, Theorem \ref{mainthm2} gives
\[
\sup_{B_{R/2}}u\leq C\left(\inf_{\mathcal{Q}_{3R}}u+ \|h_{5,\rho}(\cdot,R)^{-1}(f+I_{(\sigma(\cdot)-\rho)}u-I_{\sigma(\cdot)}u)\|_{L^\infty(B_{2\sqrt{n}R})}\right),
\]
where $h_{5,\rho}$ is defined by \eqref{h5}, replacing $\sigma(x)$ to $\sigma(x)-\rho$.
By letting $\rho\to0$, we obtain the assertion.
\end{proof}

\section{Preliminaries}\label{sec2}

\subsection{Notations}
Set
\begin{equation*}
B_r(x):=\{y\in \R^n: |y-x|<r\},
\end{equation*}
which is a open ball with center $x$, radius $r$,
and
\begin{equation*}
\mathcal{Q}_l(x):=\{y\in\R^n:|y_i-x_i|l/2\},
\end{equation*}
which is a open cube with center $x$, side length $l$.
For simplicity, $B_r$ and $\mathcal{Q}_l$ mean $B_r(0)$ and $\mathcal{Q}_l(0)$ respectively.
For $t>0$ and any open cube $\mathcal{Q}=\mathcal{Q}_l(x)$, we also write $t\mathcal{Q}:=\mathcal{Q}_{tl}(x)$.
For any measurable set $A\subset\R^n$, we denote by $|A|$ the Lebesgue measure of $A$. 
We denote $u^+:=\max\{u,0\}$ and $u^-:=\min\{-u,0\}$.

\subsection{Viscosity solutions}
Let us recall the definitions of viscosity (sub-, super-) solutions.
Here, we say that $\phi$ {\it touches} $u$ by above at $x$ whenever
\begin{equation}
u(x)=\phi(x)\quad\mbox{and}\quad u(y)\leq \phi(y)\quad\mbox{for }y\in N,\label{eq21}
\end{equation}
where $N$ is a neighborhood around $x$.
\begin{dfn}
We say that $u\in USC(B_R)\cap L^\infty(\R^n)$ is a {\it viscosity subsolution} of \eqref{eq} 
if whenever $\phi$ touches $u$ by above at $x\in B_R$ for $\phi \in C^{1,1}(\ol{N})$ and $N$ as in \eqref{eq21},   
$$
v:=
\left\{
\begin{split}
&\phi\quad {\rm in}\ N\\
&u\quad {\rm in}\ \R^n\setminus N
\end{split}
\right.
$$
satisfies that $Iv(x)\geq f(x)$.
On the other hand, $u\in LSC(B_R)\cap L^\infty(\R^n)$ is a {\it viscosity supersolution} of \eqref{eq} if $w=-u$ is a viscosity subsolution of $-Iw=-f$ in $B_R$.
Finally, $u$ is a {\it viscosity solution} of \eqref{eq} if it is both a viscosity subsolution and a viscosity supersolution of \eqref{eq}.
\end{dfn}


The {\it maximal} and {\it minimal operator} are defined by
\begin{align}
&\M^+_\nu u(x):=\int_{\R^n}(\Lambda\delta(u,x,y)^+-\lambda\delta(u,x,y)^-)\nu(x,|y|)dy,\label{max}\\
&\M^-_\nu u(x) :=\int_{\R^n}(\lambda\delta(u,x,y)^+-\Lambda\delta(u,x,y)^-)\nu(x,|y|)dy.\label{min}
\end{align}
We also have the analogous definitions of the viscosity (sub-, super-) solution of $\M^+_lu=f$ and $\M^-_lu=f$.
\begin{rem}\label{rem22}
The condition \eqref{UE} implies
\begin{equation*}
\M^-_\nu u(x)\leq Iu(x)\leq \M^+_\nu u(x).
\end{equation*}
if $u\in C^2(B_R)\cap L^\infty(\R^n)$. 
We can easily verify that if $u$ is a viscosity solution of \eqref{eq}, then $u$ is a viscosity subsolution of $\M^+_\nu u=-|f|$ and a viscosity supersolution of $\M^-_\nu u=|f|$ in $B_R$.
\end{rem}

Next, we present an analogue of Lemma 3.3 in \cite{CS09}.
\begin{prop}\label{touch}
If we have a viscosity subsolution of $\M^+_\nu u=f$ in $B_R$ and $\phi\in C^{1,1}(B_R)$ touches from above at $x\in B_R$, then $\M^+_\nu u(x)$ is defined in classical sense and $\M^+_\nu u(x)\geq f(x)$. 
\end{prop}
\begin{proof}
Let us define the auxiliary function
\[
v_r(x)=
\left\{
\begin{split}
&\phi(x)\quad {\rm in}\ B_r\\
&u(x)\quad {\rm in}\ \R^n\setminus B_r
\end{split}
\right.
\]
for any $0<r<R$. Since $\delta(v_{r_0},x,y)^+\geq\delta(v_r,x,y)^+$ for $r_0>r>0$, It is deduced from the dominate convergence theorem that
\[
\lim_{r\to0}\int_{\R^n}\delta(v_r,x,y)^+\nu(x,|y|)dy=\int_{\R^n}\delta(u,x,y)^+\nu(x,|y|)dy
\]
and so $\delta(u,x,y)^+\nu(x,|y|)$ is integrable. On the other hand, since $\delta(v_r,x,y)^-$ is increasing as $r$ decreases,
\[
\lim_{r\to0}\int_{\R^n}\delta(v_r,x,y)^-\nu(x,|y|)dy=\int_{\R^n}\delta(u,x,y)^-\nu(x,|y|)dy
\]
by the monotone convergence theorem. Hence, we obtain
\begin{align*}
\M^+_\nu u(x)=\lim_{r\to0}\M^+_\nu v_r(x)\geq f(x).
\end{align*}
The inequality $\M^+_\nu v_r(x)\geq f(x)$ implies
\[
\lambda\int_{\R^n}\delta(v_r,x,y)^-\nu(x,|y|)dy\leq\Lambda\int_{\R^n}\delta(v_r,x,y)^+\nu(x,|y|)dy-f(x).
\]
As $r\to0$, it follows that
\begin{align*}
\lambda\int_{\R^n}\delta(u,x,y)^-\nu(x,|y|)dy
\leq\lambda\int_{\R^n}\delta(u,x,y)^+\nu(x,|y|)dy-f(x)
<\infty.
\end{align*}
Thus we conclude $\delta(u,x,y)\nu(x,|y|)$ is integrable and $\M^+_\nu u(x)$ is classically defined.
\end{proof}

\subsection{Translation and scale invariance }
Let us observe an important feature of nonlocal operators related to translations and scaling.
Here we slightly modify the argument in~\cite{KL2}.
Let us consider
\begin{align}
\tilde{u}(x)&:=u(Rx+x_0)\quad\mbox{and}\label{eq2u}\\
\tilde{\nu}(x,r)&:=\frac{\nu(Rx+x_0,Rr)}{h(Rx+x_0,R)}\label{eq2nu}
\end{align}
for $x,x_0\in\R^n$, $R,r>0$.
Note that the maximal operator $\M^+_\nu u$ itself is not translation and scale invariant.
However, the translated and scaled operator $M^+_{\tilde{\nu}} u$ inherits the conditions of $\M^+_\nu u$.
\begin{prop}\label{pre2}
If $u$ is a viscosity subsolution of $\M^+_\nu u=f$ in $B_R(x_0)$, then $\tilde{u}$ is a viscosity subsolution of 
\[
\M^+_{\tilde{\nu}}\tilde{u}(x)=\frac{f(Rx+x_0)}{h(Rx+x_0,R)}\quad\mbox{in }B_1.
\]
Moreover, if $\nu$ satisfies the conditions \eqref{A1} (resp., \eqref{A2}, \eqref{A3} and \eqref{B}), then so does $\tilde{\nu}$ with the same constants.
\end{prop}

\subsection{Propseties of $h$}
Next, we recall some properties of $h(x,r)$ provided in \cite{GS17,GS18}.
To this end, we need to introduce an auxiliary function defined by
\begin{equation}\label{K}
K(x,r):=r^{-2}\int_{B_r}|y|^2dx.
\end{equation}
\begin{prop}[Lemma 2.1 and 2.2 in \cite{GS17}, Lemma 5.1 in \cite{GS18}]\label{pre3}
Suppose \eqref{A1} and \eqref{A2} hold for a function $\nu$.
Then, the following properties hold:
\begin{enumerate}
\item For any $x\in\R^n$, $K(x,\cdot)$ and $h(x,\cdot)$ are continuous and 
\[
\lim_{r\to\infty}K(x,r)=\lim_{r\to\infty}h(x,r)=0.
\]

\item $r^2K(x,r)$ and $r^2h(x,r)$ are non-decreasing for every fixed $x\in\R^n$.

\item $h(x,\cdot)$ is strictly decreasing for every fixed $x\in\R^n$ and
\[
r^{-n}K(x,r)\leq s^{-n}K(x,s)
\]
for any $x\in\R^n$ and $0<s\leq r$.

\item $\nu(x,r)\leq (n+2)|\partial B_1|^{-1}r^{-n}K(x,r).$

\item For any $0<s\leq r$,
\[
h(x,r)-h(x,s)=-\int_s^r\frac{2K(x,s)}{s}ds.
\]
\end{enumerate}
\end{prop}

\begin{prop}[Lemma 2.3 in \cite{GS17}, Lemma 5.3 in \cite{GS18}]\label{pre4}
Let $\nu$ satisfy \eqref{A1} and \eqref{A2}. The following are equivalent.
\begin{enumerate}
\item $\nu$ satisfies \eqref{A3} for some constants $\alpha,c_U>0$.


\item There is a constant $C_1>0$ such that for any $x\in\R^n$ and $r>0$, 
\[
h(x,r)\leq C_1K(x,r).
\]

\end{enumerate}

\end{prop}



\section{Alexandroff-Bakelman-Pucci Maximum Principle}

In this section, we present the Aleksandorv-Bakelman-Pucci (ABP for short) maximum principle of a viscosity subsolution $u\in USC(B_{2\sqrt{n}})\cap L^\infty(\R^n)$ of
\begin{equation}\label{subP}
\left\{
\begin{split}
&\M^+_\nu u=-f\quad\mbox{in }B_{2\sqrt{n}},\\
&u=0\hspace{15mm}\mbox{in }\R^n\setminus B_{2\sqrt{n}},
\end{split}
\right.
\end{equation}
where $\nu:\R^n\times(0,\infty)\to [0,\infty)$ is a function satisfying \eqref{A1} and \eqref{A2} (but not \eqref{A3} nor \eqref{B}) in this section.
Here, we say that $u$ is a viscosity subsolution of \eqref{subP} if $u$ is a viscosity subsolution of $\M^+_\nu u=-f$ in $B_{2\sqrt{n}}$ and satisfies $u\leq0$ in $\R^n\setminus B_{2\sqrt{n}}$.

Let us recall some notions
\begin{dfn}\label{def31}
The {\it concave envelope} $\Gamma$ in $B_{6\sqrt{n}}$ is defined by
\[
\Gamma(x):=\left\{
\begin{array}{ll}
\min\{p(x)|\ p:\mbox{affine function},p\geq u^+\mbox{ in }B_{6\sqrt{n}}\}&\mbox{in }B_{6\sqrt{n}},\\
0&\mbox{in }\R^n\setminus B_{6\sqrt{n}}.\\
\end{array}
\right.
\]
The {\it contact set} is defined as $\{u=\Gamma\}:=\{x\in B_{2\sqrt{n}};u(x)=\Gamma(x)\}.$
\end{dfn}
Now, we present a key lemma, which is an analogue of Lemma 8.1 in \cite{CS09}.
\begin{lem}\label{lem}
Let $r_k=2^{-k}/(32\sqrt{n})$ for $k=0,1,...$.
Let $\nabla\Gamma(x)$ be any element of the superdifferential of $\Gamma$ at $x$.
Then, for any $\epsilon_0\in(0,1)$, there exists $C_2=C_2(n,\lambda,\epsilon)>0$ such that if $\nu$ satisfies \eqref{A1} and \eqref{A2}, $u$ is a viscosity subsolution of \eqref{subP} and $x\in\{u=\Gamma\}$, there is a $k$ satisfying 
\begin{align}
&\quad|\{y\in B_{r_k}(x)\setminus B_{r_{k+1}}(x) : u(y)<u(x) + (y-x)\cdot\nabla\Gamma(x)-C_4f(x)r_k^2\}|\nonumber\\
&\leq\epsilon_0|B_{r_k}(x)\setminus B_{r_{k+1}}(x)|.\label{Ak}
\end{align}
\end{lem}
\begin{proof}
We first remark that we can assume $u\geq0$ in $\R^n$ without loss of generality.
If not, we consider $u^+$ instead of $u$.
From a standard argument of viscosity subsolutions, we can see $u^+$ is still a viscosity subsolution of \eqref{subP} that we omit.

For any $x \in \{u=\Gamma\},$ $u$ is touched by a plane $p$ from above at $x$.
Hence $\M^+_\nu u(x)$ is defined classically, according to lemma \ref{touch}. 
Let us observe $\delta(u,x,y)\leq0$ for $y\in\R^n$. If both $x+y\in B_{6\sqrt{n}}$ and $x-y\in B_{6\sqrt{n}}$ hold, then we derive $\delta(u,x,y)\leq\delta(p,x,y)=0$. On the other hand if either $x+y\notin B_{6\sqrt{n}}$ or $x-y\notin B_{6\sqrt{n}}$ hold, then it follows that $x+y\notin B_{2\sqrt{n}}$ and $x-y\notin B_{2\sqrt{n}}$, which imply $u(x+y)\leq0$ and $u(x-y)\leq0$, and hence $\delta(u,x,y)\leq0$.
Note that we have $\M_\nu^+u(x)\leq 0$ and so $f(x)\geq0$.

We first suppose $f(x)=0$. Then
\begin{align*}
0= \M_\nu^+ u(x)&=\int_{\R^n}\lambda\delta(u,x,y)\nu(x,|y|)dy.
\end{align*}
Hence
\[
0=\delta(u,x,y)\leq u(x+y)-u(x)-y\cdot\nabla\Gamma(x),
\]
which proves \eqref{Ak} holds for any $k=0,1,...$.

Next, we assume $f(x)>0$ in what follow.
Set
\[
A_k:=\{y\in B_{r_k}\setminus B_{r_{k+1}} : u(x+y)<u(x)+ y\cdot\nabla\Gamma(x)-C_4f(x)r_k^2\}.
\]
Suppose that \eqref{Ak} does not hold for $k=0,1,...$, in other words we suppose
\begin{equation}\label{eq31}
|A_k|\geq \epsilon_0|B_{r_k}(x)\setminus B_{r_{k+1}}(x)|\quad\mbox{for }k=0,1,....
\end{equation}
We aim to show that leads a contradiction with large $C_2$.
For $y\in A_k$, we have
\begin{align}
\delta(u,x,y)&=(u(x+y)-u(x))+(u(x-y)-u(x))\nonumber\\
&\leq(y\cdot\nabla\Gamma(x)-Mr_k^2)+(-y\cdot\nabla\Gamma(x))=-C_2f(x)r_k^2.\label{eq32}
\end{align}
Hence by using \eqref{eq31} and \eqref{eq32}, we calculate that
\begin{align*}
\int_{B_{r_0}}\delta(u,x,y)\nu(x,|y|)dy
&\leq\sum_{k=0}^\infty\int_{A_k}\delta(u,x,y)\nu(x,|y|)dy\\
&\leq-C_2f(x)\sum_{k=0}^\infty r_{k}^{2}\nu(x,r_{k})|A_k(x)|\\
&\leq -C_2f(x)\epsilon_0\sum_{k=0}^\infty r_{k}^{2}\nu(x,r_k)|B_{r_k}(x)\setminus B_{r_{k+1}}(x)|.
\end{align*}
Let $c_6=c_6(n)>0$ be a constant such that
\begin{align*}
r_{k}^2|B_{r_k}(x)\setminus B_{r_{k+1}}(x)|
&=c_6\int_{B_{r_{k-1}}(x)\setminus B_{r_{k}}(x)}|y|^2dy.
\end{align*}
Then, we have
\begin{align}
-\sum_{k=0}^\infty r_{k}^{2}\nu(x,r_k)|B_{r_k}(x)\setminus B_{r_{k+1}}(x)|\nonumber
&= -c_6\sum_{k=0}^{\infty}\int_{B_{r_{k-1}}\setminus B_{r_{k}}}|y|^2\nu(x,r_k)dy\nonumber\\
&\leq -c_6\sum_{k=0}^{\infty}\int_{B_{r_{k-1}}\setminus B_{r_{k}}}|y|^2\nu(x,|y|)dy\nonumber\\
&\leq -c_6\sum_{k=1}^{\infty}\int_{B_{r_{k-1}}\setminus B_{r_{k}}}|y|^2\nu(x,|y|)dy\nonumber\\
&=-c_6 \int_{B_{r_{0}}}|y|^2\nu(x,|y|)dy=-c_6r_0^2K(x,r_0),\nonumber
\end{align}
where the first inequality follows from
$\nu(x,|y|)\leq \nu(x,r_k)$ for $y\in B_{r_{k-1}}\setminus B_{r_{k}}$ using \eqref{A2}.
Hence we arrive at
\begin{align}
\int_{B_{r_0}}\delta(u,x,y)\nu(x,|y|)dy
&\leq-c_6r_0^2C_2f(x)\epsilon_0 K(x,r_0)\nonumber\\
&\leq -c_6(4\sqrt{n})^{-n}r_0^{2+n}C_2f(x)\epsilon_0 K(x,4\sqrt{n}),\label{eq33}
\end{align}
where we applied (3) of Proposition \ref{pre3} to the last inequality.

On the other hand, since $A_0\neq\o$ from \eqref{eq31} and $u\geq0$ in $\R^n$, fixing $z\in A_0$, it follows that
\[
-2u(x)\leq\delta(u,x,z)\leq-C_2f(x)r_0,
\]
where we applied \eqref{eq32} to the last inequality.
Hence noting $u\equiv0$ in $\R^n\setminus B_{2\sqrt{n}}$, we have
\begin{align}
\int_{\R^n\setminus B_{r_0}}\delta(u,x,y)\nu(x,|y|)dy
&\leq\int_{\R^n\setminus B_{4\sqrt{n}}}\delta(u,x,y)\nu(x,|y|)dy\nonumber\\
&=\int_{\R^n\setminus B_{4\sqrt{n}}}-2u(x)\nu(x,|y|)dy\nonumber\\
&\leq-C_2f(x)r_0^2\int_{\R^n\setminus B_{4\sqrt{n}}}\nu(x,|y|)dy.\label{eq34}
\end{align}

Consequently the estimates \eqref{eq33} and \eqref{eq34} together yield that
\begin{align*}
-f(x)
&\leq\M_\nu^+ u(x)\\
&=\int_{\R^n}\lambda\delta(u,x,y)\nu(x,|y|)dy\\
&\leq-\lambda C_2f(x)\left(c_6(4\sqrt{n})^{-n}r_0^{2+n}\epsilon_0 K(x,4\sqrt{n})+r_0^2\int_{\R^n\setminus B_2}\nu(x,|y|)dy\right)\\
&\leq-\lambda C_2f(x)\min\{c_6(4\sqrt{n})^{-n}r_0^{2+n}\epsilon_0,r_0^2\}h(x,4\sqrt{n})\\
&\leq-\lambda C_2f(x)\min\{c_6(4\sqrt{n})^{-n}r_0^{2+n}\epsilon_0,r_0^2\}(4\sqrt{n})^{-2},
\end{align*}
where the last inequality follows from \eqref{A2} and (2) of Proposition \ref{pre3}.
Hence by choosing $C_2:=2(\lambda\min\{c_6(4\sqrt{n})^{-n}r_0^{2+n}\epsilon_0,r_0^2\}(4\sqrt{n})^{-2})^{-1}$, we have a contradiction.
\end{proof}


As a consequence of Lemma \ref{lem}, we obtain the following corollary.
\begin{cor}\label{ball}
There are constants $\epsilon_1=\epsilon_1(n)$ and $C_3=C_3(n,\lambda)$ such that for any $\nu$ satisfying \eqref{A1} and \eqref{A2}, subsolution $u$ of \eqref{subP} and $x\in\{u=\Gamma\}$, there exists an $r\in (0,1/(32\sqrt{n}))$ such that
\begin{align}
&\left|\left\{y\in B_r(x)\setminus B_{r/2}(x):u(y)<u(x)+(y-x)\cdot\nabla\Gamma(x)-C_3f(x)r^2\right\}\right|\nonumber\\
\leq&\epsilon_1|B_r(x)\setminus B_{r/2}(x)|
\end{align}
and
\begin{equation}
|\nabla\Gamma(B_{r/4}(x))|\leq C_3f(x)^n|B_{r/4}(x)|.
\end{equation}
\end{cor}
\begin{proof}
Once we obtain \eqref{Ak}, the thesis follows from the same way in \cite{CS09}. See lemma 8.4 and Corollary 8.5 in \cite{CS09}.
\end{proof}

The next theorem presents the contact set $\{u=\Gamma\}$ is covered by dyadic cubes with some special properties concerned in Corollary \ref{ball}. This is proved by the same way in Theorem 8.1 \cite{CS09} that we omit here.
\begin{thm}\label{ABP}
There exist $\mu_1=\mu_1(n)$ and $C_4=C_4(n,\lambda)$ with the following property. Let $\nu$ be a function satisfying \eqref{A1} and \eqref{A2}, and $u$ be a subsolution of \eqref{subP}. There is a finite family $\{\mathcal{Q}_j\}_{j=1}^m$ of open cubes, $m\in\N$ with diameters $d_j$ such that the following hold:
\begin{align*}
&\mbox{(i)}\mathcal{Q}_j\ (j=1,...,m)\mbox{ are disjoint},\\
&\mbox{(ii)}\{u=\Gamma\}\subset \cup_{j=1}^m\overline{\mathcal{Q}}_j,\\
&\mbox{(iii)}\{u=\Gamma\}\cap\overline{\mathcal{Q}}_j\neq\o\quad(j=1,...,m),\\
&\mbox{(iv)}d_j\leq 1/(32\sqrt{n})\quad(j=1,...,m),\\
&\mbox{(v)}|\nabla\Gamma(\overline{\mathcal{Q}}_j)|\leq C_4\left(\max_{x\in\overline{\mathcal{Q}}_j}f(x)\right)^n|\mathcal{Q}_j|,\\
&\mbox{(vi)}\left|\left\{y\in 8\sqrt{n}\mathcal{Q}_j:u(y)>\Gamma(y)-C_4\left(\max_ {x\in\overline{\mathcal{Q}}_j}f(x)\right)d_j^2\right\}\right|\geq\mu_1|\mathcal{Q}_j|.
\end{align*}
\end{thm}


As a consequence of Theorem \ref{ABP}, we obtain an upper bound of viscosity subsolutions which is a nonlocal version of ABP maximum principle (see Theorem 3.2 in \cite{CafCab} for the second order version of ABP maximum principles).
\begin{cor}\label{abp}
There is a constant $C_5=C_5(n,\lambda)$ such that for any function $\nu$ satisfying \eqref{A1} and \eqref{A2}, and subsolution $u$ of \eqref{subP}, we have
\[
\sup_{B_{2\sqrt{n}}}u\leq C_5 \left(\sum_{j=1}^m(\max_{x\in\overline{\mathcal{Q}}_j}f(x))^n|\mathcal{Q}_j|\right)^{1/n}
\]
where $\{\mathcal{Q}_j\}_{j=1}^m$ is as in theorem \ref{ABP}.
\end{cor}
\begin{proof}
Without loss of generality, we can assume $\sup_{B_1}u>0$ and $u(x_0)=\sup_{B_1}u$ for $x_0\in B_{2\sqrt{n}}$.
We write $M_3:=u(x_0)$.
We first show that 
\begin{equation*}
B_{M_3/(12\sqrt{n})}\subset \nabla\Gamma(\{u=\Gamma\})
\end{equation*}
holds.
To see this, fix any $p\in B_{M_3/(12\sqrt{n})}$ and consider $x_1\in \overline{B}_{6\sqrt{n}}$ such that
\[
u(x_1)^++\langle p,x_1\rangle=\max_{B_{6\sqrt{n}}}\{u(x)^++\langle p,x\rangle\}.
\]
Then, it follows that
\begin{align*}
u(x_0)&\leq u(x_1)+\langle p,x_1-x_0\rangle\\
&\leq u(x_1)+(12\sqrt{n})|p|\\
&<u(x_1)+M_3=u(x_1)+u(x_0),
\end{align*}
and hence $u(x_1)>0$.
This implies $x_1\in B_{2\sqrt{n}}$ and the affine function $u(x_1)+\langle p,x-x_1\rangle$ touches $u$ from above at $x=x_1$.
Hence from Definition \ref{def31}, we have $x_1\in\{u=\Gamma\}$ and $p\in\nabla\Gamma(x_1)\subset\nabla\Gamma(\{u=\Gamma\})$,
which prove the claim.
Moreover, that yields that
\begin{equation}
\left(\sup_{B_{2\sqrt{n}}}u\right)^n=c_7|B_{M_3/(12\sqrt{n})}|\leq c_7|\nabla\Gamma(\{u=\Gamma\}).\label{eq351}
\end{equation}
for a constant $c_7>0$ depending only on $n$.

On the other hand, from (vi) of theorem \ref{ABP}, we arrive at
\begin{align*}
|\nabla\Gamma(\{u=\Gamma\})|
&\leq\sum_{j=1}^m|\nabla\Gamma(\overline{\mathcal{Q}}_j)|\\
&\leq C_4\left(\max_{x\in\overline{\mathcal{Q}}_j}
f(x)\right)^n|\mathcal{Q}_j|
\end{align*}
which together with \eqref{eq351}, proves the assertion of the corollary with $C_5=(c_7C_4)^{1/n}$.
\end{proof}



\section{Barrier function}

In this section, we introduce a barrier function which is used in the proof of the weak Harnack inequality in sect \ref{sect5}.
Note that the proofs of Lemma 9.1 of \cite{CS09}, Lemma 3.5 of \cite{KKL16} and Lemma 4.4 of \cite{KL} cannot be directly applied to our equations because these are based on the fact that $|x|^{-p}$ for $p\geq n$ is not locally integrable. 
To deal with this gap, we developed a new approach based on the properties of $h$.

\begin{lem}\label{lem2}
Assume that $\nu$ satisfies (A).
There is a constant $p=p(n,\lambda,\Lambda,c_U,\alpha)>0$ such that for any $s_1>0$ the function
\[
f(x)=\min( (s_1/2)^{-p},|x|^{-p})
\]
satisfies that
\begin{equation}\label{bar}
\M^-_\nu f(x)\geq0\quad\mbox{for }|x|>s_1.
\end{equation}
\end{lem}
\begin{proof}
By rotational symmetry, it is enough to check \eqref{bar} holds for $x=re_1=(r,0,...,0)$ with $r>s_1$.

Set
\begin{align*}
\tau_0&:=\sqrt{\frac{\lambda}{2n\Lambda(p+4)}},\\
I_1&:=\int_{B_{\tau_0 r}}(\lambda\delta(f,x,y)^+-\Lambda\delta(f,x,y)^-)\nu(x,|y|)dy\quad\mbox{and}\\
I_2&:=\int_{\R^n\setminus B_{\tau_0 r}}(\lambda\delta(f,x,y)^+-\Lambda\delta(f,x,y)^-)\nu(x,|y|)dy.
\end{align*}
By Taylor's theorem, we have the following elementary estimates for any $a>b>0$ and $q>0$:
\begin{align*}
(a+b)^{-q}&\geq a^{-q}\left(1-q\frac{b}{a}\right),\\
(a+b)^{-q}+(a-b)^{-q}&\geq2a^{-q}+q(q+1)b^2a^{-q-2}.
\end{align*}
For $y\in B_{r/2}$, we compute
\begin{align}
\delta(f,x,y)&=|x+y|^{-p}+|x-y|^{-p}-2|x|^{-p}\nonumber\\
&=(r^2+|y|^2+2ry_1)^{-p/2}+(r^2+|y|^2-2ry_1)^{-p/2}-2r^{-p}\nonumber\\
&\geq2(r^2+|y|^2)^{-p/2}+p(p+2)\frac{y_1^2}{r^2}(r^2+|y|^2)^{-p/2-2}-2r^{-p}\nonumber\\
&\geq-pr^{-p}\left(\frac{|y|^2}{r^2}+(p+2)\frac{y_1^2}{r^2}-\frac{1}{2}(p+2)(p+4)\frac{y_1^2|y|^2}{r^4}\right)\nonumber.
\end{align}
Hence, we have
\begin{align*}
&\quad\lambda\delta(f,x,y)^+-\Lambda\delta(f,x,y)^-\\
&\geq pr^{-p}\left(\lambda(p+2)\frac{y_1^2}{r^2}-\Lambda\frac{|y|^2}{r^2}-\frac{\Lambda}{2}(p+2)(p+4)\frac{y_1^2|y|^2}{r^4}\right)
\end{align*}
for $y\in B_{r/2}$.
Choosing $p$ large enough such that
\begin{equation*}
\frac{\lambda(p+2)}{2n}-\Lambda>0,
\end{equation*}
we have
\begin{align}
I_1
&\geq pr^{-p}\int_{B_{\tau_0 r}}\left(\lambda(p+2)\frac{y_1^2}{r^2}-\Lambda\frac{|y|^2}{r^2}-\frac{\Lambda}{2}(p+2)(p+4)\frac{y_1^2|y|^2}{r^4}\right)\nu(x,|y|)dy\nonumber\\
&=pr^{-p}\int_{B_{\tau_0 r}}\left\{\frac{\lambda(p+2)}{n}-\Lambda-\frac{\Lambda(p+2)(p+4)}{2n}\cdot\frac{|y|^2}{r^2}\right\}\frac{|y|^2}{r^2}\nu(x,|y|)dy\nonumber\\
&\geq pr^{-p}\int_{B_{\tau_0 r}}\left\{\frac{\lambda(p+2)}{2n}-\frac{\Lambda(p+2)(p+4)}{2n}\cdot\tau_0^2\right\}\frac{|y|^2}{r^2}\nu(x,|y|)dy\nonumber\\
&= pr^{-p}\frac{\lambda (p+2)}{4n}\int_{B_{\tau_0 r}}\frac{|y|^2}{r^2}\nu(x,|y|)dy\nonumber\\
&= \frac{\lambda p(p+2)\tau_0^2}{4nr^p}K(x,\tau_0 r)\nonumber\\
&\geq\frac{\lambda p(p+2)\tau_0^2}{4nr^p}C_1^{-1}h(x,\tau_0 r)\nonumber\\
&\geq c_8pr^{-p}h(x,\tau_0 r),\nonumber
\end{align}
where we applied (2) of Proposition \ref{pre4} to the third inequality and the fact that $(p+2)/(p+4)$ is bounded below to the last inequality so that $c_8>$ depends only on $n,\lambda,\Lambda,c_U$ and $\alpha$.

On the other hand, since $\delta(f,x,y)\geq -2f(x)=-2r^{-p}$, we have
\begin{align*}
I_2
&\geq\int_{\R^n\setminus B_{\tau_0 r}}-\Lambda\delta(f,x,y)^-\nu(x,|y|)dy\\
&\geq -2\Lambda r^{-p}\int_{\R^n\setminus B_{\tau_0 r}}\nu(x,|y|)dy
\geq -2\Lambda r^{-p}h(x,\tau_0 r).
\end{align*}

Combining the lower bounds of $I_1$ and $I_2$, we obtain
\begin{align*}
\M^-_\nu f(x)
&\geq (c_8p-2\Lambda)r^{-p}h(x,\tau_0 r).
\end{align*}
Letting $p>0$ large enough, \eqref{bar} is proved.
\end{proof}


\begin{cor}\label{bar2}
There is a function $\Phi$ such that for any $\nu$ satisfying (A),
\begin{align*}
&\mbox{(i) }\Phi(x)=0\mbox{ in }\R^n\setminus B_{2\sqrt{n}},\\
&\mbox{(ii) }\Phi(x)\geq2\mbox{ in }\mathcal{Q}_{3}\quad\mbox{and}\\
&\mbox{(iii) }\M^-_\nu\Phi(x)\geq -\psi(x)\mbox{ in }\R^n,
\end{align*}
where $\psi$ is some positive function supported in $\overline{B}_{1/4}$.
\end{cor}

\begin{proof}
Let $p$ as in lemma \ref{lem2}.
We construct a function $\Phi$ such that
\begin{equation*}
\Phi(x)
=c_4\left\{
\begin{array}{ll}
q(x)&\mbox{in }B_{1/8},\\
|x|^{-p}-(2\sqrt{n})^{-p}&\mbox{in }B_{2\sqrt{n}}\setminus B_{1/8},\\
0&\mbox{in }\R^n\setminus B_{2\sqrt{n}},
\end{array}
\right.
\end{equation*}
where we choose $c_9:=2/\{(3\sqrt{n}/2)^{-p}-(2\sqrt{n})^{-p}\}$ and
\[
q(x)=-\frac{p8^{p}}{2}|x|^2+\left\{\left(p+1\right)8^{p}-(2\sqrt{n})^{-p}\right\}
\]
so that $\Phi$ is $C^{1,1}$ on $\partial B_{ 1/8}$ and $\Phi\geq2$ in $B_{3\sqrt{n}/2}$.
Since $\mathcal{Q}_3\subset B_{3\sqrt{n}/2}$, (ii) holds.

Now it remains to check (iii).
Let $f$ be as in lemma \ref{lem2} with $s_0=1/4$.
Then, for any $x\in B_{2\sqrt{n}}\setminus B_{1/8}$, $f$ touches $\Phi$ from below at $x$, and so we have
$\M^-_\nu\Phi(x)\geq \M^-_\nu f(x)\geq0$.
In the case $x\in\R^n\setminus B_{2\sqrt{n}}$, $\delta(\Phi,x,y)\geq0$ is clear by the construction of $\Phi$, which implies $M_\nu \Phi(x)\geq0$.
\end{proof}



\section{Point to Measure Estimate}\label{sect5}

We present in this section, the point to measure estimate which is a version of the weak Harnack inequality and important in the proof of the H\"older estimate and the Harnack inequality.  


\begin{lem}\label{pe1}
There exist constants $\epsilon_2>0$, $\mu_2\in(0,1)$ and $M_4>1$ depending only on $n,\lambda,\Lambda,c_U$ and $\alpha$ such that if (A) is satisfied and $u\in C(B_1)\cap L^\infty(\R^n)$ is a function such that
\begin{align*}
&\mbox{(i) }u\geq0\quad\mbox{in }\R^n,\\
&\mbox{(ii) }\inf_{\mathcal{Q}_{3}}u\leq1\quad\mbox{and}\\
&\mbox{(iii)}u\mbox{ is a viscosity supersolution of}\\
&\quad \M^-_\nu u(x)\leq \epsilon_2\quad\mbox{in }B_{2\sqrt{n}},
\end{align*}
then, it follows that
\begin{equation}\label{PE1}
|\{u\leq M_4\}\cap\mathcal{Q}_1|>\mu_2|\mathcal{Q}_1|.
\end{equation}
\end{lem} 

\begin{proof}
We consider $v=\Phi-u$ with the function $\Phi$ as in corollary \ref{bar2}.
Then we easily verify that $v$ is a viscosity subsolution of
\[
\M^+_\nu v(x)=  -\epsilon_2+\M^-_\nu\Phi(x)\quad\mbox{in }B_{2\sqrt{n}}
\]
and $v\leq0$ in $\R^n\setminus B_{2\sqrt{n}}$.
We apply theorem \ref{ABP} and Corollary \ref{abp} to derive
\[
\sup_{B_{2\sqrt{n}}}v
\leq \frac{C_5}{2\sqrt{n}}\left(\sum_{j=1}^m\left\{\max_{\overline{\mathcal{Q}}_j}\epsilon_2-\M^-_{\tilde{\nu}}\Phi(x)\right\}^n|\mathcal{\mathcal{Q}}_j|\right)^{1/n},
\]
where $\{\mathcal{Q}_j\}_{j=1}^m$ is a family of disjoint open cubes satisfying
\begin{equation}
\left|\left\{y\in 8\sqrt{n}\mathcal{Q}_j:v(y)>\Gamma(y)-C_4\left\{\max_{\overline{\tilde{\mathcal{Q}}_j}}\epsilon_2-\M^-_{\tilde{\nu}}\Phi(2\sqrt{n}x)\right\}d_j^2\right\}\right|\geq\mu_1|\mathcal{Q}_j|.
\label{pe11}
\end{equation}
for $j=1,...,m$, $\mathrm{diam}(\mathcal{Q}_j)=:d_j\in (0,1/(32\sqrt{n}))$ and the concave envelope $\Gamma$ of $v$ defined in Definition \ref{def31}.
By using (iii) of Corollary \ref{bar2}, we calculate as follow:
\begin{align*}
\sup_{B_{2\sqrt{n}}}v
&\leq C_5\left(\sum_{j=1}^m\left\{\max_{\overline{\mathcal{Q}}_j}\epsilon_2-\M^-_{\tilde{\nu}}\Phi(x)\right\}^n|\mathcal{\mathcal{Q}}_j|\right)^{1/n}\\
&\leq C_5\left(\sum_{j=1}^m\left\{\max_{\overline{\mathcal{Q}}_j}\epsilon_2+\psi(x)\right\}^n|\mathcal{Q}_j|\right)^{1/n}\\
&\leq C_5
\left\{\epsilon_2|B_{2\sqrt{n}+1/(32\sqrt{n})}|^{1/n}
+\|\psi\|_\infty\left (\sum_{j=1}^m(\max_{\overline{\mathcal{Q}}_j}\chi_{B_{1/4}}(x))^n|\mathcal{Q}_j|\right)^{1/n}
\right\}\\
&\leq c_{10}\epsilon_2+c_{10}\left (\sum_{\mathcal{Q}_j\cap B_{1/4}\neq\o}|\mathcal{Q}_j|\right)^{1/n},
\end{align*}
where $c_{10}:=C_5\max\{|B_{2\sqrt{n}+1/(32\sqrt{n})}|^{1/n},\|\psi\|_\infty\}$.
Since $\inf_{\mathcal{Q}_{3}}u\leq 1$ and $\inf_{\mathcal{Q}_{3}}\Phi\geq 2$, we have $\sup_{B_{2\sqrt{n}}}v\geq1$.
Hence by choosing $\epsilon_2=(2c_{10})^{-1}$, we conclde that
\begin{equation}\label{below1}
\sum_{j=1}^{m'}|\mathcal{Q}_j|\geq \frac{1}{(2c_{10})^n},
\end{equation}
where we have relabeled $\mathcal{Q}_j$ so that they intersect with $B_{1/4}$ for $j=1,...,m'$ and not for $j=m'+1,...,m$.

Let $x^j$ be a center of cubes $\mathcal{Q}_j$.  Recall that $d_j\leq 1/(32\sqrt{n})$. We have $|x^j|\leq 1/4+1/(64\sqrt{n})$ and $|y|\leq |x^j|+1/8\leq 1/2$ for $y\in8\sqrt{n}\mathcal{Q}_j$. Therefore $\mathcal{Q}_j\cap B_{1/4}\neq\o$ implies $8\sqrt{n}\mathcal{Q}_j\subset B_{1/2}$. 

From \eqref{pe11}, we have
\begin{equation}\label{below2}
|\{y\in8\sqrt{n}\mathcal{Q}_j:v(y)\geq-M_5|\geq \mu_1|\mathcal{Q}_j|
\end{equation}
for some constant $M_5=M_5(n,\lambda,\Lambda,c_U,\alpha)$ since $\Gamma$ is non-negative and we estimate
\begin{align*}
C_4\left\{\max_{\overline{\mathcal{Q}}_j}\epsilon_2-\M^-_{\tilde{\nu}}\Phi(x)\right\}d_j^2
&\leq C_4\left\{\max_{\overline{\mathcal{Q}}_j}\epsilon_2+\psi(x)\right\}d_j^2\\
&\leq \frac{C_4(\epsilon_2+\|\psi\|_\infty)}{(32\sqrt{n})^2}=:M_5,
\end{align*}
where we applied (iii) of Corollary \ref{bar2} to the first inequality and $d_j\in(0,1/(32\sqrt{n}))$ to the second inequality.

By applying the Besicovitch covering theorem to the family of the cubes $8\sqrt{n}\mathcal{Q}_j$ (see for instance Theorem 18.1c in \cite{DiBe16}),
there exists subfamily $8\sqrt{n}\mathcal{Q}_{jk}$ which covers $\{x^j\}_{j=1}^{m'}$ and at most $4^n$ cubes overlap in $\R^n$.
Suppose $x^j$ is included in $8\sqrt{n}\mathcal{Q}_{jk}$.
Then, at least one-$2^n$ th part of $\mathcal{Q}_j$ has to be also covered by $\{8\sqrt{n}\mathcal{Q}_{jk}\} $ (see Figure \ref{Fig1}).
\begin{figure}[h]
\begin{center}
\includegraphics[keepaspectratio,scale=0.3]{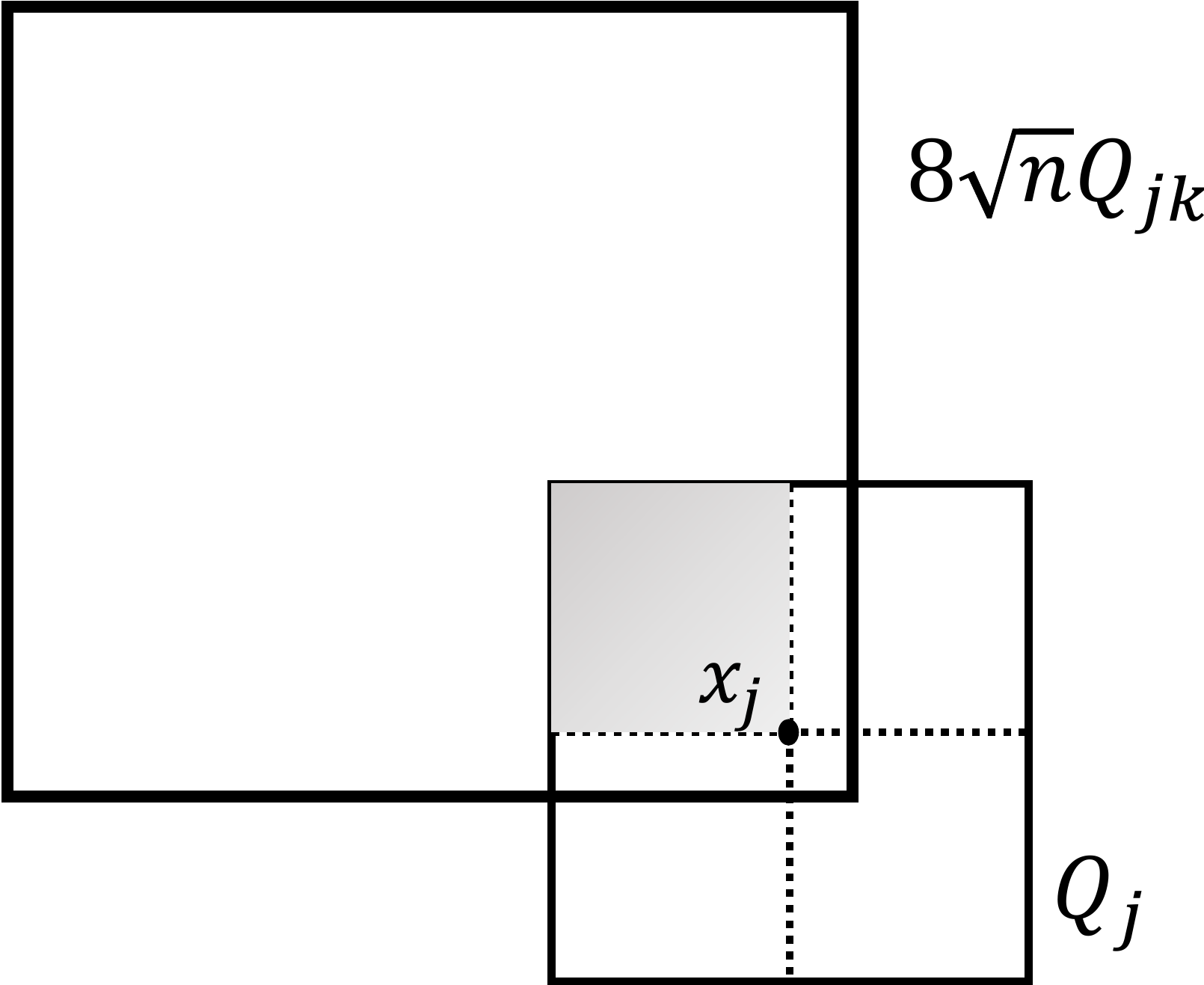}
\caption{one-$2^n$ th part of $\mathcal{Q}_j$ covered by $\mathcal{Q}_{jk}$}\label{Fig1}
\end{center}
\end{figure}

\noindent
Hence, we have
\begin{equation}\label{below3}
\sum_{k}|8\sqrt{n}\mathcal{Q}_{jk}|\geq 2^{-n}\sum_{j=1}^{m'}|\mathcal{Q}_j|.
\end{equation}

We put $\cup_{s=1}^{4^n}\{8\sqrt{n}\mathcal{Q}_{jk}^s\}_k=\{8\sqrt{n}\mathcal{Q}_{jk}\}_k$ such that for each $s$, $\{8\sqrt{n}\mathcal{Q}^s_{jk}\}_k$ is disjoint.
Since $8\sqrt{n}\mathcal{Q}_{jk}^s$ is contained in $B_{1/2}$, it follows that
\begin{align*}
|\{y\in B_{1/2}:v(y)\geq-M_5\}|
&\geq|\{y\in\cup_{k}8\sqrt{n}\mathcal{Q}_{jk} : v(y)\geq-M_5\}|\\
&\geq 4^{-n}\sum_{s=1}^{4^n}|\{y\in\cup_{k}8\sqrt{n}\mathcal{Q}_{jk}^s: v(y)\geq-M_5\}|\\
&= 4^{-n}\sum_{s,k}|\{y\in8\sqrt{n}\mathcal{Q}_{jk}^s: v(y)\geq-M_5\}|.
\end{align*}
Combine \eqref{below1}, \eqref{below2} and \eqref{below3}, to obtaine
\begin{align*}
|\{y\in B_{1/2}:v(y)\geq-M_5\}|
\geq 4^{-n}\mu_1\sum_{s,k}|8\sqrt{n}\mathcal{Q}^s_{jk}|
&\geq 8^{-n}\mu_1\sum_{j=1}^{m'}|\mathcal{Q}_j|\\
&\geq \mu_1(16c_{10})^{-n}.
\end{align*}
Note that $v(y)\geq -M_5$ implies 
\begin{align*}
u(y)
\leq \Phi(y)+M_5\leq\|\Phi\|_{L^\infty(\R^n)}+M_5.
\end{align*}
Set $M_4:=\|\Phi\|_{L^\infty(\R^n)}+M_5$ and $\mu_2:=\mu_1(16c_{10})^{-n}$. Then, because $B_{1/2}\subset\mathcal{Q}_1$, we have
\begin{equation*}
|\{y\in\mathcal{Q}_1:u(y)\leq M_5\}|> \mu_2|\mathcal{Q}_1|.
\end{equation*}
\end{proof}


\begin{lem}\label{pe2}
There exists constants $C_6,\rho_0>0$ depending only on $n,\lambda,\Lambda,c_U$ and $\alpha$ such that for any $u$ satisfying the hypothesis of lemma \ref{pe1},
it follows that
\begin{equation}\label{PE2}
|\{u> t\}\cap\mathcal{Q}_1|\leq C_6t^{-\rho_0}|\mathcal{Q}_1|\quad\mbox{for }t>0.
\end{equation}
\end{lem}

\begin{proof}
It is enough to prove
\begin{equation}\label{PE3}
|\{u>M_4^k\}\cap\mathcal{Q}_1|\leq (1-\mu_2)^k|\mathcal{Q}_1|
\end{equation}
for $k=1,2,...$, where $M_4$ and $\mu_2$ are as in lemma \ref{pe1}. \eqref{PE2} follows from that with taking $C_6:=(1-\mu_2)^{-1}$ and $\rho_0>0$ such that $1-\mu_2=M_4^{-\rho_0}$.

For $k=1$, \eqref{PE3} is just \eqref{PE1}. Suppose that \eqref{PE3} holds for $k-1$ and let
\begin{equation*}
A:=\{u>M_4^k\}\cap\mathcal{Q}_1,\quad B:=\{u>M_4^{k-1}\}\cap\mathcal{Q}_1.
\end{equation*}
That will be proved if we show that
\begin{equation}\label{cube2}
|A|\leq(1-\mu_2)|B|,
\end{equation}
by applying the Calder\'on-Zygmund cube decomposition (see lemma 4.2 in \cite{CafCab}).
We need to check the following condition:
if $\mathcal{Q}$ is a dyadic cube such that
\begin{equation}\label{cube}
|A\cap\mathcal{Q}|\geq(1-\mu_2)|\mathcal{Q}|,
\end{equation}
then $\tilde{\mathcal{Q}}\subset B$ for the predecessor $\tilde{\mathcal{Q}}$ of $\mathcal{Q}$ i.e. $\mathcal{Q}$ is one of the $2^n$ cubes obtained from dividing $\tilde{\mathcal{Q}}$.
On the contrary, let us suppose that there exists $x\in\mathcal{Q}_1$ and $R\in(0,1)$ such that \eqref{cube} holds for $\mathcal{Q}_R(x_0)$ and $\tilde{\mathcal{Q}}\nsubseteq B$ also holds for the predecessor $\tilde{Q}$ of $\mathcal{Q}_R(x_0)$.
We show that leads to a contradiction.

Let us consider the function
\[
\tilde{u}(x):=\frac{u(Rx+x_0)}{M^{k-1}_4}
\]
and $\tilde{\nu}(x,r)$ as in \eqref{eq2nu}.
Clearly $\tilde{u}\geq0$.
We can also check that $\inf_{\mathcal{Q}_{3}}\tilde{u}\leq1$ from $\tilde{\mathcal{Q}}\nsubseteq B$
and $\tilde{u}$ is a viscosity supersolution of 
\begin{align*}
\M^-_{\tilde{\nu}}\tilde{u}(x)= \frac{\epsilon_2}{M_4^{k-1}} \quad\mbox{in }B_{2\sqrt{n}}
\end{align*}
according to Proposition \ref{pre2}.
Since $\epsilon_2/M_4^{k-1}\leq\epsilon_2$, $\tilde{u}$ is under the hypothesis of lemma \ref{pe1}, and so we have
\begin{equation*}
\mu_2|\mathcal{Q}_1|<|\{\tilde{u}\leq M_4\}\cap\mathcal{Q}_{1}|=R^{-n}|\{u\leq M_4^k\}\cap \mathcal{Q}_R(x_0)|.
\end{equation*}
Therefore $|\mathcal{Q}_R(x_0)\setminus A|>\mu_2|\mathcal{Q}_R(x_0)|$, which contradicts \eqref{cube}.
Hence $A$ and $B$ satisfies the hypothesis of Lemma 4.2 in \cite{CafCab} and we obtain \eqref{cube2}.
\end{proof}


\begin{thm}\label{pe4}
Let $u\geq0$ in $\R^n$ and $\M^-_\nu u\leq |f|$ in $B_{2\sqrt{n}R}$. Then, it follows that
\begin{equation}\label{PE5}
|\{u> t\}\cap \mathcal{Q}_R|\leq C_6\left(\inf_{\mathcal{Q}_{3R}}u+\|h(\cdot,R)^{-1}f(\cdot)\|_{L^\infty(B_{2\sqrt{n}R})}\right)^{\rho_0}t^{-\rho_0}|\mathcal{Q}_R|\quad\mbox{for }t>0.
\end{equation}
where $C_6>0$ and $\rho_0>0$ are from Lemma \ref{pe2}.
\end{thm}

\begin{proof}
We consider for any $\eta_0>0$,
\begin{align*}
u_{\eta_0}(x)
&:=\frac{u(Rx)}{\inf_{\mathcal{Q}_{3R}}u+\eta_0+ \|h(\cdot,R)^{-1}f\|_{L^\infty(B_{2\sqrt{n}R})}/\epsilon_2}\quad\mbox{and}\\
f_{\eta_0}(x)
&:=\frac{h(Rx,R)^{-1}|f(Rx)|}{\inf_{\mathcal{Q}_{3R}}u+\eta_0+ \|h(\cdot,R)^{-1}f(\cdot)\|_{L^\infty(B_{2\sqrt{n}R})}/\epsilon_2}
\end{align*}
where $\epsilon_2$ is as in lemma \ref{pe1}. Let $\tilde{\nu}$ be as in \eqref{eq2nu}.
Then, we have $\inf_{\mathcal{Q}_3}u_{\eta_0}\leq 1$ and $|f_{\eta_0}|\leq\epsilon_2$ and according to Proposition \ref{pre2}, $u_{\eta_0}$ is a viscosity supersolution of
\begin{align*}
\M^-_{\tilde{\nu}} u_{\eta_0}(x)=f_{\eta_0} \quad\mbox{in }\mathcal{Q}_{2\sqrt{n}}.
\end{align*}
Thus $u_{\eta_0}$ is under the hypothesis of Lemma \ref{pe2} and so
\begin{equation*}
 R^n|\{u_{\eta_0}> t\}\cap \mathcal{Q}_{1}|\leq C_6t^{-\rho}|\mathcal{Q}_R|.
\end{equation*}
By taking $s=t/(\inf_{\mathcal{Q}_{3R}}u+\eta_0+ \|h(\cdot,R)^{-1}f(\cdot)\|_{L^\infty(B_{2\sqrt{n}R})}/\epsilon_2)$, it follows that
\begin{equation*}
|\{u> s\}\cap \mathcal{Q}_{R}|\leq C_6\left(\inf_{\mathcal{Q}_{3R}}u+\eta_0+\|h(\cdot,R)^{-1}f(\cdot)\|_{L^\infty(B_{2\sqrt{n}R})}/\epsilon_2\right)^{\rho_0}s^{-\rho_0}|\mathcal{Q}_R|
\end{equation*}
Letting $\eta_0\to0$, \eqref{PE5} is proved.
\end{proof}

Theorem \ref{pe4} implies the following weak Harnack inequality:
\begin{thm}
For any $\rho_1\in(0,\rho_0)$, there exists $C_{7}>0$, depending only on $n,\lambda,\Lambda,c_U$, $\alpha$ and $\rho_1$ with the following property. Let $u\geq0$ in $\R^n$ and $\M^-_\nu u\leq |f|$ in $B_{2\sqrt{n}R}$. Then, it follows that
\begin{equation}\label{PE5}
\|u\|_{L^{\tilde{\rho}}(\mathcal{Q}_R)}\leq C_{10}\left(\inf_{\mathcal{Q}_{3R}}u+\|h(\cdot,R)^{-1}f(\cdot)\|_{L^\infty(B_{2\sqrt{n}R})}\right).
\end{equation}
\end{thm}
\begin{proof}
The thesis immediately follows from the identity
\[
\|u\|_{L^{\rho_1}(\mathcal{Q}_R)}^{\rho_1}=\rho_1\int_{0}^{\infty}t^{\rho_1-1}|\{u> t\}\cap \mathcal{Q}_R|dt.
\]
and the decay estimate of the distribution function of $u$ in Theorem \ref{pe4}.
\end{proof}


\section{H\"older estimate}
Now, we are ready to prove Theorem \ref{main thm1}.
Noting Remark \ref{rem22}, it suffice to show the following lemma:
\begin{lem}\label{lem61}
There exist $C_7>0$ and $\mu_1>0$ such that if the function $\nu$ satisfies (A) and and $u$ is a viscosity subsolution of $\M^+_\nu u(x)=-|f|$ in $B_R$ and a viscosity supersolution of $\M^-_\nu u(x)=|f|$ in $B_R$ for $R>0$, then it follows that
\[
\| u\|_{C^{\mu_1}(B_{R/2})}\leq C_7R^{-{\mu_1}}(\|u\|_{L^\infty(\R^n)}+\|h(\cdot,R)^{-1}f(\cdot)\|_{L^\infty(B_R)}).
\]
\end{lem}

\begin{proof}
Fix $\epsilon_3>0$ such that
\begin{equation}
0<c_{11}\epsilon_3\leq\frac{1}{2(2C_6)^{1/\rho_0}},\label{eq64}
\end{equation}
where we defined $c_{11}>0$ by
\begin{equation}
c_{11}:=2c_U(1+4(\alpha+1)\alpha^{-2}\Lambda (n+2)),\label{eq68}
\end{equation}
and $\mu_1$ such that
\begin{align}
0<\mu_1&\leq \alpha,\label{eq62}\\
\frac{2^{\mu_1}}{\alpha-\mu_1}-\frac{1}{\alpha}&\leq \epsilon_3\quad\mbox{and}\label{eq63}\\
0<2(1-(2\sqrt{n})^{-\mu_1})&\leq \frac{1}{2(2C_6)^{1/\rho_0}}.\label{eq65}
\end{align}
Set
\begin{equation}
\theta_1:=2(1-(2\sqrt{n})^{-\mu_0}).\label{eq67}
\end{equation}

Without loss of generality, we can assume $R=1$, $\|f\|_\infty\leq \epsilon_3$ and $\|u\|_\infty\leq1/2$.
Indeed, in general case, we consider
\begin{align*}
\hat{u}(x)&:=\frac{u(Rx)}{2\|u\|_{L^\infty(\R^n)}+\|h(\cdot,R)^{-1}f(\cdot)\|_{L^\infty(B_R)}/\epsilon_3},\\
\hat{f}(x)&:=\frac{h(Rx,R)^{-1}f(Rx)}{2\|u\|_{L^\infty(\R^n)}+\|h(\cdot,R)^{-1}f(\cdot)\|_{L^\infty(B_R)}/\epsilon_3}
\end{align*}
and $\tilde{\nu}$ defined by \eqref{eq2nu}.
Then, according to Proposition \ref{pre2}, we can observe that $\|\hat{f}\|_\infty\leq \epsilon_3$, $\|\hat{u}\|_\infty\leq1/2$ and $\hat{u}$ satisfies the hypothesis of Lemma \ref{lem61} for $f=\hat{f}$, $\nu=\tilde{\nu}$ and $R=1$.
It is enough to prove the assertion for $\hat{u}$, which is $\|\hat{u}\|_{C^{\eta_1}(B_{1/2})}\leq C$ for some constant $C>0$ depending only on $n,\lambda, \Lambda, c_U$ and $\alpha$,
since this implies
\begin{align*}
\| u\|_{C^{\eta_1}(B_{R/2})}
&=\|\hat{u}\|_{C^{\eta_1}(B_{1/2})}R^{-{\eta_1}}(2\|u\|_{L^\infty(\R^n)}+\|h(\cdot,R)^{-1}f(\cdot)\|_{L^\infty(B_R)}/\epsilon_3)\\
&\leq C\max\left\{2,\frac{1}{\epsilon_3}\right\}R^{-{\eta_1}}(\|u\|_{L^\infty(\R^n)}+\|h(\cdot,R)^{-1}f(\cdot)\|_{L^\infty(B_R)})
\end{align*}
and we complete the proof of the assertion for $u$.

For any fixed $x_0\in B_{1/2}$, we aim to show that there exist two sequences $\{m_k\}_{k=1}^\infty$ and $\{M_k\}_{k=1}^\infty$ such that
\begin{equation}\label{eq61}
m_k\leq u\leq M_k\quad\mbox{in }\mathcal{Q}_{(2\sqrt{n})^{-k}}(x_0)\quad\mbox{and}\quad M_k-m_k=(2\sqrt{n})^{-\mu_1 k}.
\end{equation}
This construction proves the assertion of the theorem as we have
\[
|u(x)-u(x_0)|\leq M_k-m_k=(2\sqrt{n})^{-\mu_1 k}\leq 2\sqrt{n}|x-x_0|^{\mu_1}.
\]
for $x\in\mathcal{Q}_{(2\sqrt{n})^{-k}}(x_0)\setminus\mathcal{Q}_{(2\sqrt{n})^{-k-1}}(x_0)$.

For $k=1$, we choose $m_0=-1/2$ and $M_0=1/2$.
Let us construct the sequences $M_k$ and $m_k$ for $k=2,3,...$ by induction.
Let $j\geq2$.
Suppose that we already have the sequences $m_k$ and $M_k$ satisfying \eqref{eq61} for $k\leq j-1$.
We need to find $m_{j}$ and $M_{j}$ which still satisfy \eqref{eq61}.

Set $R_j:=1/(3(2\sqrt{n})^{j})$. It is easy to see that one of these two inequalities holds:
\begin{align*}
\left|\left\{u\geq \frac{M_{j-1}+m_{j-1}}{2}\right\}\cap \mathcal{Q}_{R_j}(x_0)\right|&\geq \frac{|\mathcal{Q}_{R_j}|}{2}\quad\mbox{and}\\
\left|\left\{u\leq \frac{M_{j-1}+m_{j-1}}{2}\right\}\cap \mathcal{Q}_{R_j}(x_0)\right|&\geq \frac{|\mathcal{Q}_{R_j}|}{2}.
\end{align*}

We first consider that the first inequality holds.
Set
\[
v(x):=\frac{u(R_jx+x_0)-m_{j-1}}{(M_{j-1}-m_{j-1})/2}
\]
which satisfies $v(x)\geq 0$ in $\mathcal{Q}_{6\sqrt{n}}$ and $|\{v\geq 1\}\cap \mathcal{Q}_{1}|\geq |\mathcal{Q}_{1}|/2$.
For $\tilde{\nu}$ as in \eqref{eq2nu} with $R=R_j$, $v$ is a viscosity supersolution of
\[
\M^-_{\tilde{\nu}} v(x)=2c_U\epsilon_3\quad\mbox{in }B_{2\sqrt{n}}
\]
since we have
\begin{align*}
\frac{\epsilon_3}{h(R_jx+x_0,R_j)(M_{j-1}-m_{j-1})/2}
&\leq 2\epsilon_3c_U (6\sqrt{n})^{-\alpha}(2\sqrt{n})^{(\mu_1-\alpha) (j-1)}\leq 2c_U\epsilon_3
\end{align*}
for $x\in B_{2\sqrt{n}}$, where we used \eqref{A3} to the first inequality and \eqref{eq62} to the last inequality.

From the inductive hypothesis, for any $i\leq j$ and $x\in\mathcal{Q}_{3(2\sqrt{n})^{i}}\setminus\mathcal{Q}_{3(2\sqrt{n})^{i-1}}$, we have
\begin{align*}
v(x)
\geq \frac{(m_{j-i}-m_{j-1})}{(M_{j-1}-m_{j-1})/2}
&1\geq\frac{(m_{j-i}-M_{j-i}+M_{j-1}-m_{j-1})}{(M_{j-1}-m_{j-1})/2}\\
&\geq 2(1-(2\sqrt{n})^{\mu_1 (i-1)})\\
&\geq 2\left(1-\left(\frac{2}{3\sqrt{n}}\right)^{\mu_1} |x|^{\mu_1}\right).
\end{align*}
Hence we conclude that 
\begin{equation}
v(x)\geq -2\left(1-\left(\frac{2}{3\sqrt{n}}\right)^{\mu_1} |x|^{\mu_1}\right)^-\quad\mbox{in }\R^n.\label{eq69}
\end{equation}

We remark that $v$ is not necessarily nonnegative in $\R^n$.
Hence we have to consider $w:=v^+$, to apply Theorem \ref{pe4}.
Because $v=w$ in $\mathcal{Q}_{6\sqrt{n}}$, we can easily verify that $w$ is a viscosity supersolution of $\M^-_{\tilde{\nu}} w(x)= 2c_U\epsilon_3+\M^+_{\tilde{\nu}}(w-v)(x)$ in $B_{2\sqrt{n}}$, where for $x\in B_{2\sqrt{n}}$, $\M^+_{\tilde{\nu}}
(w-v)(x)$ is defined in classical sense and
\begin{align*}
\M^+_{\tilde{\nu}}(w-v)(x)
&=\int_{\R^n}2\Lambda v(x+y)^-\tilde{\nu}(x,|y|)dy\\
&=\int_{\R^n\setminus B_{\sqrt{n}}}2\Lambda v(x+y)^-\tilde{\nu}(x,|y|)dy.
\end{align*}
Setting
\[
\tilde{h}(x,r):=\int_{\R^n}\left(1\wedge\frac{|y|^2}{r^2}\right)\tilde{\nu}(x,|y|)dy,
\]
we compute
\begin{align*}
\M^+_{\tilde{\nu}}(w-v)(x)
&\leq\int_{\R^n\setminus B_{\sqrt{n}}}4\Lambda\left(1-\left(\frac{2}{3\sqrt{n}}\right)^{\mu_1}|x+y|^{\mu_1}\right)^-\tilde{\nu}(x,|y|)dy\\
&\leq 4\Lambda (n+2)|\partial B_1|^{-1}\int_{\R^n\setminus B_{\sqrt{n}}}\left(1-\left(\frac{2}{3\sqrt{n}}\right)^{\mu_1}(|x|+|y|)^{\mu_1}\right)^-\frac{\tilde{h}(x,|y|)}{|y|^n}dy\\
&\leq 4\Lambda (n+2)|\partial B_1|^{-1}\int_{\R^n\setminus B_{\sqrt{n}}}\left(1-\left(\frac{2}{\sqrt{n}}\right)^{\mu_1}|y|^{\mu_1}\right)^-\frac{c_U\tilde{h}(x,1)}{|y|^{n+\alpha}}dy\\
&=4c_U\Lambda (n+2)\int_{\sqrt{n}}^\infty\left(\left(\frac{2}{\sqrt{n}}\right)^{\mu_1} r^{-1-\alpha+\mu_1}-r^{-1-\alpha}\right)dr\\
&=4n^{-\alpha/2}c_U\Lambda (n+2)\left(\frac{2^{\mu_1}}{\alpha-\mu_1}-\frac{1}{\alpha}\right)\\
&\leq4c_U\Lambda (n+2)\frac{2(\alpha+1)}{\alpha^2}\epsilon_3,
\end{align*}
where we have applied \eqref{eq69} to the first inequality, (4) of Proposition \ref{pre3} and $K(x,|y|)\leq h(x,|y|)$ to the second inequality, \eqref{A3} to the third inequality and \eqref{eq63} to the last inequality.
Consequently, $w$ is a viscosity supersolution of $M_{\tilde{\nu}} w(x)= c_{11}\epsilon_3$ in $B_{2\sqrt{n}}$ for $c_{11}$ as in \eqref{eq68}.

Applying Theorem \ref{pe4} to $w$, we have
\begin{align*}
C_6\left(\inf_{\mathcal{Q}_3}w+c_{11}\epsilon_3\right)^{\rho_0}|\mathcal{Q}_1|&\geq |\{w\geq1\}\cap \mathcal{Q}_1|\geq\frac{|\mathcal{Q}_{1}|}{2}.
\end{align*}
Combining this inequality, \eqref{eq64}, \eqref{eq65} and \eqref{eq67}, we have $\inf_{\mathcal{Q}_3}w\geq1/(2(2C_6)^{1/{\rho_0}})\geq \theta_1$.
Thus letting $M_{j}:=M_{j-1}$ and $m_{j}:=m_{j-1}+\theta_1(M_{j-1}-m_{j-1})/2$, we have $m_{j}\leq u\leq M_{j}$ in $\mathcal{Q}_{(2\sqrt{n})^{-j}}(x_0)$.
Moreover, it follows from \eqref{eq67} that $M_{j}-m_{j}=(1-\theta_1/2)(2\sqrt{n})^{-\mu_1 (j-1)}=6^{-\mu_1 j}$.

On the other hand, if $|\{u\leq (M_{j-1}+m_{j-1})/2\}\cap \mathcal{Q}_{R_j}(x_0)|\geq |\mathcal{Q}_{R_j}|/2$, we define
\[
v(x):=\frac{M_{j-1}-u(R_jx+x_0)}{(M_{j-1}-m_{j-1})/2}
\]
and repeat the previous argument with some minor modifications.
\end{proof}


\section{Harnack inequality}
From the discussion in Remark \ref{rem22}, Theorem \ref{mainthm2} is reduced to the following lemma:
\begin{lem}\label{he}
There exists $C_8>0$ with the following properties.
Let (A) and (B) hold.
If $u$ is non-negative in $\R^n$, a viscosity subsolution of $\M^+_\nu u=|f|$ and a viscosity supersolution of $\M^-_\nu u=-|f|$ in $B_{2\sqrt{n}R}$ for $R>0$,
then it follows that
\begin{equation}\label{HI}
\sup_{B_{R/2}}u\leq C_8\left(\inf_{\mathcal{Q}_{3R}}u+ \|h(\cdot,R)^{-1}f\|_{L^\infty(B_{2\sqrt{n}R})}\right).
\end{equation}
\end{lem}

\begin{proof}
First we note that it suffice to prove our assertion \eqref{HI} in the simplified case $R=1$, $\inf_{\mathcal{Q}_3}u\leq 1$ and $|f|\leq1$ in $B_{2\sqrt{n}}$.
In fact, in the general case, we may consider
\begin{align*}
\hat{u}(x)&:=\frac{u(Rx)}{\inf_{\mathcal{Q}_{3R}}u+ \|h(\cdot,R)^{-1}f\|_{L^\infty(B_{2\sqrt{n}R})}},\\
\hat{f}(x)&:=\frac{h(Rx,R)^{-1}f(Rx)}{\inf_{\mathcal{Q}_{3R}}u+ \|h(\cdot,R)^{-1}f\|_{L^\infty(B_{2\sqrt{n}R})}}
\end{align*}
and $\tilde{\nu}$ defined by \eqref{eq2nu}.
From Proposition \ref{pre2}, we can observe that $\hat{u}$ satisfies the hypothesis of Lemma \ref{he} for $f=\hat{f}$, $\nu=\tilde{\nu}$ and $R=1$.
Moreover $\inf_{\mathcal{Q}_{3}}\hat{u}\leq1$ and $|\hat{f}|\leq1$ hold.
Hence we restrict ourselves to $R=1$, $\inf_{\mathcal{Q}_3}u\leq 1$ and $|f|\leq1$ in what follow.

We will use the family $\{V_i\}_{i=1}^{N_0}$ of open cones later such that $N_0$ depends only on the dimension $n$, $\bigcup_{i=1}^{N_0}\overline{V}_i=\R^n$, the vertex of $V_i$ is $0$ and the angle of $V_i$ is $\pi/3$ i.e. For the axis $z_i\in\partial B_1$ of $V_i$, $y\in V_i$ holds if and only if
\[
\left\langle\frac{y}{|y|},z_i \right\rangle\geq \mathrm{cos} \frac{\pi}{6}=\frac{\sqrt{3}}{2}
\]
(see figure \ref{Fig3}).
\begin{figure}[h]
\begin{center}
\includegraphics[scale=0.3]{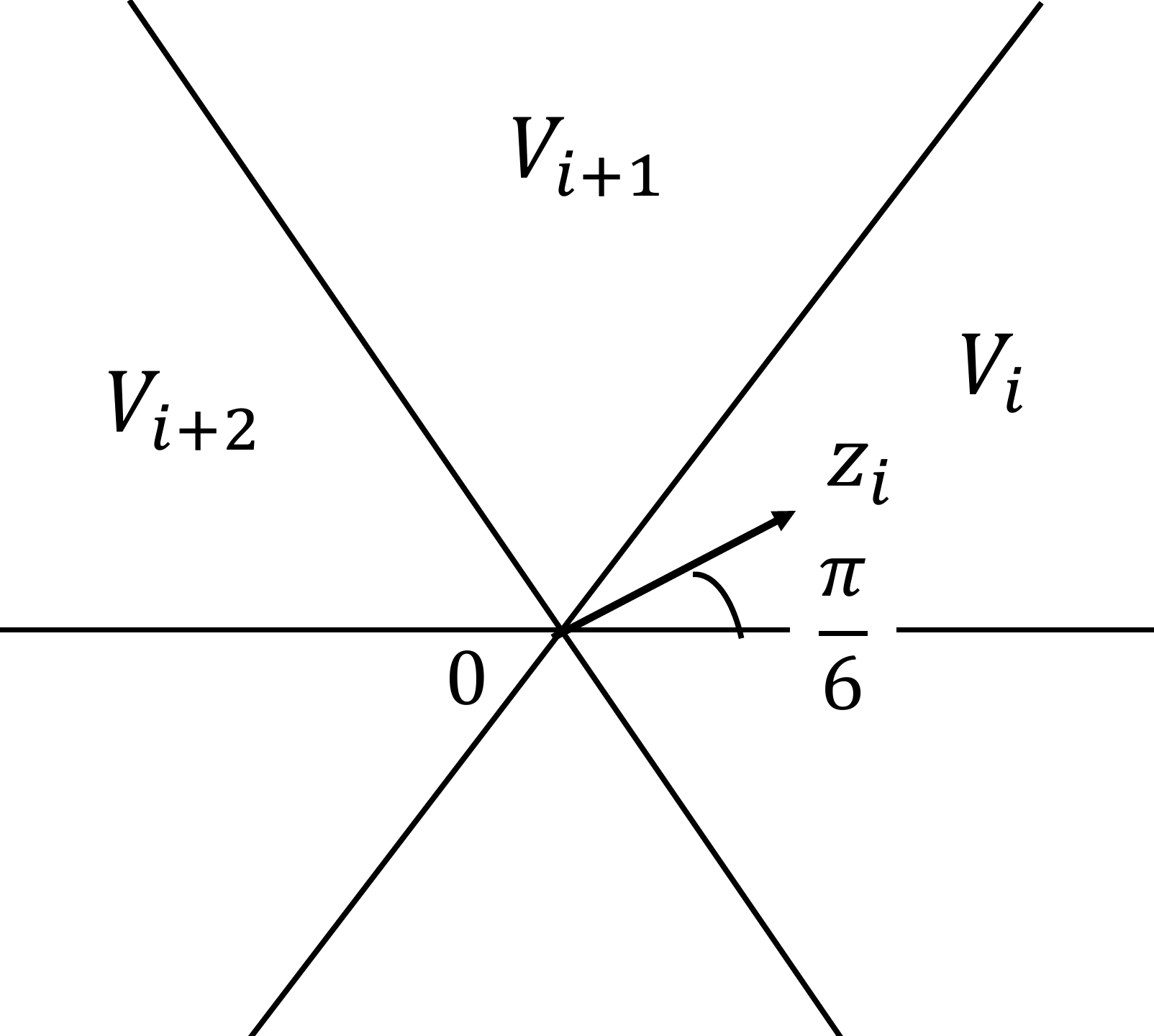}
\caption{$V_i$ and $z_i$}\label{Fig3}
\end{center}
\end{figure}
Next, let us fix some constants which depend only on $n,\lambda,\Lambda,c_U,\alpha$ and $c^*$.
Let $\gamma_1:=n/\rho_0$ and $c_{12}$ be the positive constants defined by
\begin{equation}
c_{12}:=1+2\Lambda c^*N_0\lambda^{-1}\left(1+\left(\Lambda128+\lambda\left(\frac{8}{3}\right)^{2} \right)2^{\gamma_1}\right),\label{he7}
\end{equation}
where $\rho_0$ is from Lemma \ref{pe2}.
We also choose $\theta_2\in(0,1)$ and $\tau_1\in(0,1/4)$ small and $t_0\geq1$ large to satisfy
\begin{align}
C_6\left(\left(1-\frac{\theta_2}{2}\right)^{-\gamma_1}-1\right)^{\rho_0}
&\leq\frac{1}{4},\label{he6}\\
C_6(c_{12}c_U\tau_1^{\alpha})^{\rho_0}2^{\rho_0}&\leq\frac{1}{4}\label{he11}\\
4^{\rho_0}C_6t_0^{-\rho_0}
&<\frac{(\tau_1\theta_2)}{2(2\sqrt{n})^n},\label{he8}
\end{align}
where $C_6$ is from Lemma \ref{pe2}.


Let $t_u$ be the smallest value of $t$ such that
\begin{equation*}
u(x)\leq h_t(x):=t(1-|x|)^{-\gamma_1}\quad\mbox{for every }x\in B_{1
}.
\end{equation*}
Then, there must be an $x_0\in B_1$ for which $u(x_0)=h_{t_u}(x_0)$.

Our aim is to show the upper bound $t_u<t_0$, which implies
\begin{equation*}
u(x)\leq t_0(1-|x|)^{-\gamma_1}\leq t_02^{\gamma_1} \quad\mbox{in }B_{1/2},
\end{equation*}
and hence the thesis follows.

On contrary, suppose that there exists $u$ such that $t_u>t_0$.
Set
\begin{align}
A&:=\left\{u>\frac{u(x_0)}{2}\right\},\quad d_0:=\frac{1-|x_0|}{2}\quad\mbox{and}\nonumber\\
\epsilon_4&:=\left(1-\frac{\theta_2}{2}\right)^\gamma\left(\frac{\theta_2 d}{8}\right)^2\frac{1}{u(x_0)}.\label{he4}
\end{align}
Then, Theorem \ref{pe4} gives that
\begin{align}
|A\cap \mathcal{Q}_1|
\leq C_6\left(\frac{2}{u(x_0)}\right)^{\rho_0}\left(\inf_{\mathcal{Q}_3}u+\|f\|_{L^\infty(B_{2\sqrt{n}})}\right)^{\rho_0}|\mathcal{Q}_1|
&\leq 4^{\rho_0}C_6t_u^{-\rho_0}d_0^{\rho_0\gamma_1}\nonumber\\
&\leq4^{\rho_0}C_6t_0^{-\rho_0}d_0^{n}\label{he1},
\end{align}
where we used $\inf_{\mathcal{Q}_3}u\leq1$, $\|f\|_{L^\infty(B_{2\sqrt{n}})}\leq1$ and $u(x_0)=t_u/(2d_0)^{\gamma_1}$ to the second inequality and $\gamma_1=n/\rho_0$ to the third inequality.
Hence to get a contradiction, it suffice to show the inequality opposite to \eqref{he1}.

For every $x\in B_{\theta_2 d_0}(x_0)\subset B_1$, we have 
\begin{equation}
u(x)\leq h_{t_u}(x)\leq t_u(d_0-\theta_2 d_0)^{-\gamma_1}=u(x_0)\left(1-\frac{\theta_2}{2}\right)^{-\gamma_1}\label{he16}
\end{equation}
Set
\begin{equation*}
v(x):=\left(1-\frac{\theta_2}{2}\right)^{-\gamma_1}u(x_0)-u(x).
\end{equation*}
which satisfies $v\geq0$ in $B_{\theta_2 d_0}(x_0)$, $\M_\nu^- v(x)\leq 1$ and $\M_\nu^+ v(x)\geq -1$ in $B_{2\sqrt{n}}$.
Since $v$ does not satisfies the hypothesis of Theorem \ref{pe4}, we also need to consider $w:=v^+$ instead.
Because of the non-negativity of $v$ in $B_{\theta_2 d_0}(x_0)$, we can easily verify that $w$ is a viscosity supersolution of
\[
\M^-_\nu w(x)=1+\M^+_\nu(w-v)(x)\quad\mbox{in }B_{\tau_1\theta_2 d_0}(x_0)
\]
where $\M^+_\nu(w-v)(x)$ is defined in classical sense and
\begin{align}
\M^+_\nu(w-v)(x)
&=\int_{\R^n}\Lambda (v(x+y)^-+v(x-y)^-)\nu(x,|y|)dy\nonumber\\
&=2\Lambda\int_{\R^n\setminus B_{(1-\tau_1)\theta_1 d_0}} v(x+y)^-\nu(x,|y|)dy\nonumber.
\end{align}
Next, we fix $x\in B_{\tau_1\theta_2 d_0}(x_0)$ and estimate $\M^+_\nu (w-v)(x)$ from above.
Consider the cones $V_i$ for $i=1,...,N_0$, which we constructed at the beginning of the proof.
Here we note that 
for any $y_1,y_2\in V_i$ we have
\begin{equation}
\left\langle\frac{y_1}{|y_1|},\frac{y_2}{|y_2|} \right\rangle
\geq \mathrm{cos} \frac{\pi}{3}=\frac{1}{2}\label{he3}
\end{equation}
and $B_{1/2}(z_i)\subset V_i$ (see Figurer \ref{Fig2}).
\begin{figure}[h]
\begin{center}
\includegraphics[scale=0.3]{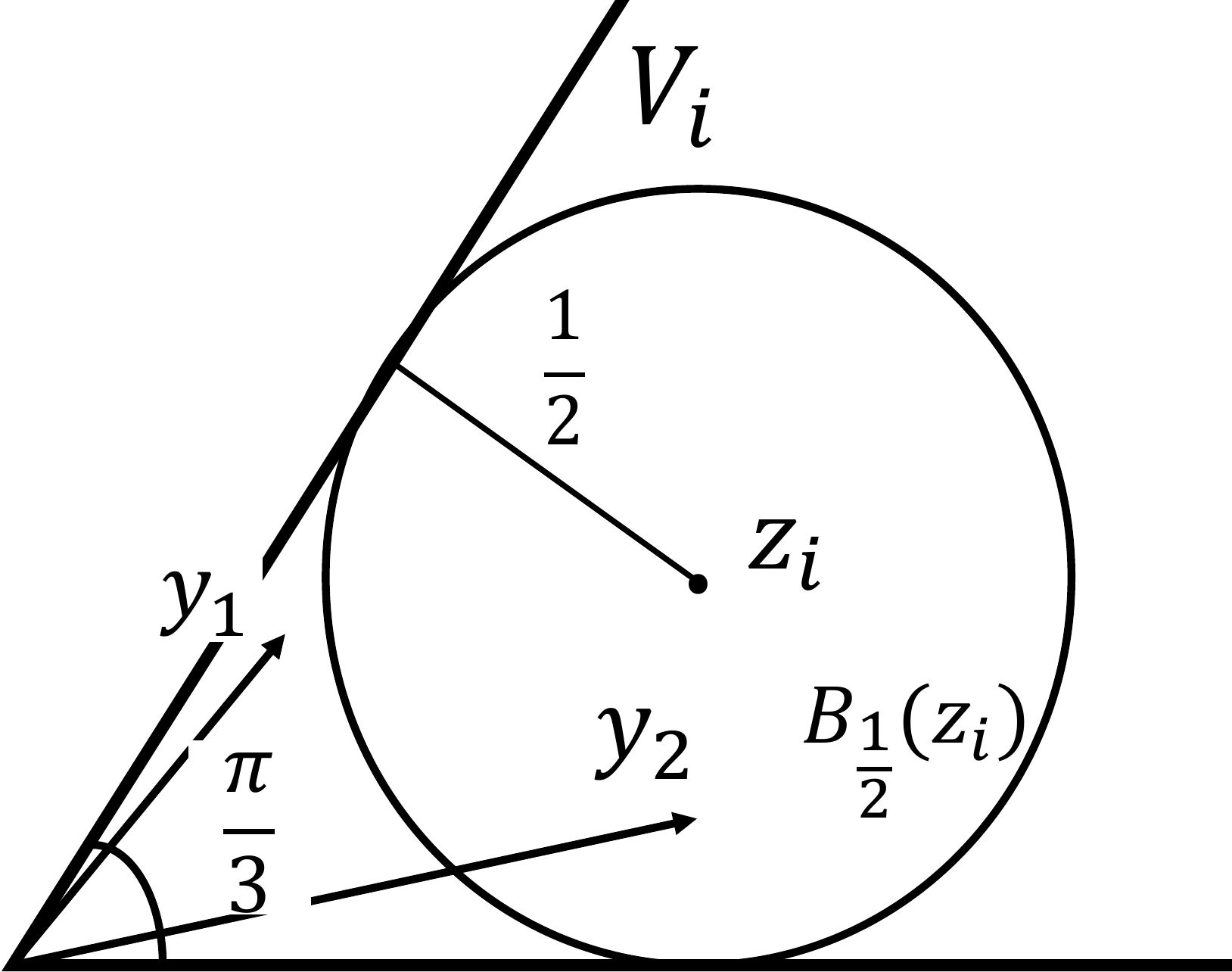}
\caption{$V_i$ and $B_{1/2}(z_i)$}\label{Fig2}
\end{center}
\end{figure}
From symmetry to translation and scaling, we still have $B_{\theta_2 d_0/8}(x_i)\subset x+V_i$ for $x_i:=x+(\theta_2 d_0/4)z_i$.
\noindent
For $y\in V_i\setminus B_{(1-\tau_1)\theta_2 d_0}$ and $\tilde{x}\in B_{\theta_2 d_0/8}(x_i)\subset x+V_i$, it follows that
\begin{equation}
|\tilde{x}-x|
\leq|\tilde{x}-x_i|+|x_i-x|\leq \left(\frac{1}{8}+\tau_1\right)\theta_2 d_0\leq |y|,\label{he12}
\end{equation}
where we used $(1/8+\tau_1)\leq (1-\tau_1)$ for $\tau_1\in(0,1/4)$, and
\begin{align}
|y-\tilde{x}+x|^2
&=|y|^2-2\langle y,\tilde{x}-x\rangle+|\tilde{x}-x|^2\nonumber\\
&\leq|y|^2-|y||\tilde{x}-x|+|\tilde{x}-x|^2\leq|y|^2,\label{he5}
\end{align}
where we applied \eqref{he3} to the second inequality and \eqref{he12} to the last inequality.
Since $u(z)+|z-x_i|^2/\epsilon_4 \to\infty$ as $|z|\to\infty$, we can consider $\tilde{x}_i\in \R^n$ such that
\[
u(\tilde{x}_i)+\frac{1}{\epsilon_4}|\tilde{x}_i-x_i|^2=\inf_{z\in\R^n}\left\{u(z)+\frac{1}{\epsilon_4}|z-x_i|^2\right\}.
\]
We can take $\tilde{x}_i\in B_{\theta_2 d_0/8}(x_i)\subset \tilde{V}_i$ since
\begin{align*}
|\tilde{x}_i-x_i|^2
&\leq \epsilon_4(u(x_i)-u(\tilde{x}_i))\\
&\leq  \epsilon_4\left(\left(1-\frac{\theta_2}{2}\right)^{-\gamma_1} u(x_0)-u(\tilde{x}_i)\right)\\
&\leq  \epsilon_4\left(1-\frac{\theta_2}{2}\right)^{-\gamma_1} u(x_0)= \left(\frac{\theta_2 d_0}{8}\right)^2,
\end{align*}
where we used \eqref{he16} to the second inequality, $u\geq0$ to the third inequality and \eqref{he4} to the last equality.
Since $-|z-x_i|^2/\epsilon_4$ touches $u$ from below at $z=\tilde{x}_i$,
Proposition \ref{touch} implies $\M^-_\nu u(\tilde{x}_1)$ is defined in classical sense and we can calculate as follow:
\begin{align*}
&\quad\int_{\R^n}\delta(u,\tilde{x}_i,y)^-\nu(\tilde{x}_i,|y|)dy\\
&\leq\int_{B_{\theta_2 d_0}}\frac{2}{\epsilon_4}|y|^2\nu(\tilde{x}_i,|y|)dy+\int_{\R^n\setminus B_{\theta_2 d_0}}(2u(\tilde{x}_i)-u(\tilde{x}_i+y)-u(\tilde{x}_i-y)^+\nu(\tilde{x}_i,|y|)dy\\
&\leq\frac{2(\theta_2 d_0)^{2}}{\epsilon_4}K(\tilde{x}_i,\theta_2 d_0)+2u(\tilde{x}_i)\int_{\R^n\setminus B_{\theta_2 d_0}}\nu(\tilde{x}_i,|y|)dy\\
&\leq 128\left(1-\frac{\theta_2}{2}\right)^{-\gamma_1} u(x_0)h(\tilde{x}_i,\theta_2 d_0),
\end{align*}
where we applied \eqref{he16} and \eqref{he4} to the last inequality.
Since $\M^-_\nu u(\tilde{x}_i)\leq 1$, we have
\begin{align}
\int_{\R^n}\delta(u,\tilde{x}_i,y)^+\nu(\tilde{x}_i,|y|)dy
&\leq\lambda^{-1}\left(1+\Lambda128\left(1-\frac{\theta_2}{2}\right)^{-\gamma_1} u(x_0)h(\tilde{x}_i,\theta_2 d_0)\right).\nonumber
\end{align}
We also have
\begin{align*}
&\quad\int_{\R^n}\delta(u,\tilde{x}_i,y)^+\nu(\tilde{x}_i,|y|)dy\\
&\geq\int_{\R^n\setminus B_{(7/8-2\tau_1)\theta_2 d_0}}\left(u(\tilde{x}_i+y)-2u(x_i)\right)^+\nu(\tilde{x}_i,|y|)dy\\
&\geq\int_{\R^n\setminus B_{(7/8-2\tau_1)\theta_2 d_0}}\left(u(\tilde{x}_i+y)-2\left(1-\frac{\theta_2}{2}\right)^{-\gamma_1} u(x_0)\right)^+\nu(\tilde{x}_i,|y|)dy\\
&\geq\int_{\R^n\setminus B_{(7/8-2\tau_1)\theta_2 d_0}}\left(v(\tilde{x}_i+y)^--\left(1-\frac{\theta_2}{2}\right)^{-\gamma_1}u(x_0)\right)\nu(\tilde{x}_i,|y|)dy\\
&\geq\int_{\R^n\setminus B_{(7/8-2\tau_1)\theta_2 d_0}}v(\tilde{x}_i+y)^-\nu(\tilde{x}_i,|y|)dy\\
&\quad-\left(1-\frac{\theta_2}{2}\right)^{-\gamma_1} u(x_0)h\left(\tilde{x}_i,\left(\frac{7}{8}-2\tau_1\right)\theta_2 d_0\right).
\end{align*}
Hence
\begin{align}
&\quad\int_{\R^n\setminus B_{(7/8-2\tau_1)\theta_2 d_0}}v(\tilde{x}_i+y)^-\nu(\tilde{x}_i,|y|)dy\nonumber\\
&\leq \lambda^{-1}\left(1+\left(\Lambda128h(\tilde{x}_i,\theta_2 d_0)+\lambda h\left(\tilde{x}_i,\left(\frac{7}{8}-2\tau_1\right)\theta_2 d_0\right)\right)\left(1-\frac{\theta_2}{2}\right)^{-\gamma_1} u(x_0)\right)\nonumber\\
&\leq \lambda^{-1}\left(1+\left(\Lambda128+\lambda\left(\frac{7}{8}-2\tau_1\right)^{-2} \right)\left(1-\frac{\theta_2}{2}\right)^{-\gamma_1}h(\tilde{x}_i,\theta_2 d_0) u(x_0)\right)\nonumber\\
&\leq \lambda^{-1}\left(1+\left(\Lambda128+\lambda\left(\frac{8}{3}\right)^{2} \right)2^{\gamma_1}\right)h(\tilde{x}_i,\theta_2 d_0)u(x_0),\label{he13}
\end{align}
where we applied (2) of Proposition \ref{pre3} to the second inequality, $h(\tilde{x}_i,\theta_2 d_0)\geq h(\tilde{x}_i,1)=1$, $u(x_0)\geq t_0\geq1$, $\tau_1\in(0,1/4)$ and $\theta_2\in(0,1)$ to the last inequality.
For $y\in V_i\setminus B_{(1-\tau_1)\theta_2 d_0}$, we also have
\begin{align}
\nu(x,|y|)\leq c^*\frac{\nu (\tilde{x}_i,|y|)h(x,\theta_2 d_0)}{h(\tilde{x}_i,\theta_2 d_0)}\leq c^*\frac{\nu(\tilde{x}_i,|y-\tilde{x}_i+x|)h(x,\theta_2 d_0)}{h(\tilde{x}_i,\theta_2 d_0)},\label{he14}
\end{align}
where we used \eqref{B} and \eqref{he12} to the first inequality and \eqref{A2} and \eqref{he5} to the second inequality.
Noting that we have $B_{(7/8-2\tau_1)\theta_2 d_0}\leq B_{(1-\tau_1)\theta_2 d_0}(x-\tilde{x}_i)$ from \eqref{he12} for each $i=1,...,N_0$, we calculate as follows:
\begin{align*}
&\quad\int_{\R^n\setminus B_{(1-\tau_1)\theta_2 d_0}} v(x+y)^-\nu(x,|y|)dy\\
&\leq\sum_{i=1}^N\int_{V_i\setminus B_{(1-\tau_1)\theta_2 d_0}}v(x+y)^-\nu(x,|y|)dy\\
&\leq c^*\sum_{i=1}^N\int_{V_i\setminus B_{(1-\tau_1)\theta_2 d_0}}v(\tilde{x}_i+y-\tilde{x}_i+x)^-\frac{\nu(\tilde{x}_i,|y-\tilde{x}_i+x|)h(x,\theta_2 d_0)}{h(\tilde{x}_i,\theta_2 d_0)}dy\\
&\leq c^*\sum_{i=1}^N\int_{\R^n\setminus B_{(1-\tau_1)\theta_2 d_0}(x-\tilde{x}_i)}v(\tilde{x}_i+y)^-\frac{\nu(\tilde{x}_i,|y-\tilde{x}_i+x|)h(x,\theta_2 d_0)}{h(\tilde{x}_i,\theta_2 d_0)}dy\\
&\leq c^*\sum_{i=1}^N\int_{\R^n\setminus B_{(7/8-2\tau_1)\theta_2 d_0}}v(\tilde{x}_i+y)^-\frac{\nu(\tilde{x}_i,|y-\tilde{x}_i+x|)h(x,\theta_2 d_0)}{h(\tilde{x}_i,\theta_2 d_0)}dy\\
&\leq c^*N\lambda^{-1}\left(1+\left(\Lambda128+\lambda\left(\frac{8}{3}\right)^{2} \right)2^{\gamma_1}\right)h(x,\theta_2 d_0)u(x_0),
\end{align*}
where we applied \eqref{he14} to the second inequality and \eqref{he13} to the last inequality.
Consequently, $w$ is a viscosity supersolution of
\[
\M_\nu^-w(x)=c_{12}h(x,\theta_2 d_0)u(x_0)\quad\mbox{in }B_{\tau_1\theta_2 d_0}(x_0)
\]
since it follows that
\begin{align*}
1+2\Lambda c^*N\lambda^{-1}\left(1+\left(\Lambda128+\lambda\left(\frac{8}{3}\right)^{2} \right)2^{\gamma_1}\right)h(x,\theta_2 d_0)u(x_0)\leq c_{12}h(x,\theta_2 d_0)u(x_0),
\end{align*}
where we applied $h(x,\theta_2 d_0)\geq h(x,1)=1$ and $u(x_0)\geq t_0\geq1$ and $c_{12}$ is from \eqref{he7}.

Applying Theorem \ref{pe4} to $w$ with $R=\tau_1\theta_2 d_0/(2\sqrt{n})$, we have
\begin{align*}
&\quad\left|\left\{u\leq \frac{u(x_0)}{2}\right\}\cap \mathcal{Q}_{\tau_1 \theta_2 d_0/(2\sqrt{n})}(x_0)\right|\\
&=\left|\left\{w\geq u(x_0)\left(\left(1-\frac{\theta_2}{2}\right)^{-\gamma_1}-\frac{1}{2}\right)\right\}\cap \mathcal{Q}_{\tau_1 \theta_2 d_0/(2\sqrt{n})}(x_0)\right|\\
&\leq C_6 \left(\frac{\tau_1\theta_2 d_0}{2\sqrt{n}}\right)^n\left(u(x_0)\left(\left(\left(1-\frac{\theta_2}{2}\right)^{-\gamma_1}-1\right)+\left\|\frac{c_{12}u(x_0)h(\cdot,\theta_2 d_0)}{h(\cdot,\tau_1\theta_2 d_0)}\right\|_{L^\infty(B_{\tau_1 \theta_2 d_0})}\right)\right)^{\rho_0}\\
&\quad\cdot\left(u(x_0)\left(\left(1-\frac{\theta_2}{2}\right)^{-\gamma_1}-\frac{1}{2}\right)\right)^{-\rho_0}\\
&\leq C_6 \left(\frac{\tau_1\theta_2 d_0}{2\sqrt{n}}\right)^n\left(\left(\left(\left(1-\frac{\theta_2}{2}\right)^{-\gamma_1}-1\right)+c_{12}\tau_1^\alpha\right)\right)^{\rho_0}2^{\rho_0}\\
&\leq C_6\left(\frac{\tau_1\theta_2 d_0}{2\sqrt{n}}\right)^n\left(\left(\left(1-\frac{\theta_2}{2}\right)^{-\gamma_1}-1\right)^{\rho_0}+(c_{12}c_U\tau_1^{\alpha})^{\rho_0}\right)2^{\rho_0},
\end{align*}
where we used \eqref{A3} to the second inequality.
From \eqref{he6}, we obtain
\[
C_6\left(\frac{\tau_1\theta_2 d_0}{2\sqrt{n}}\right)^n\left(\left(1-\frac{\theta_2}{2}\right)^{-\gamma_1}-1\right)^{\rho_0}2^{\rho_0}\leq\frac{1}{4}|\mathcal{Q}_{\tau_1 \theta_2 d_0/(2\sqrt{n})}(x_0)|
\]
and from \eqref{he11},
\begin{align*}
C_6\left(\frac{\tau_1\theta_2 d_0}{2\sqrt{n}}\right)^n(c_{12}c_U\tau_1^{\alpha})^{\rho_0}2^{\rho_0}\leq \frac{1}{4}|\mathcal{Q}_{\tau_1 \theta_2 d_0/(2\sqrt{n})}(x_0)|.
\end{align*}
Therefore,
\[
\left|\left\{u\leq \frac{u(x_0)}{2}\right\}\cap \mathcal{Q}_{\tau_1 \theta_2 d_0/(2\sqrt{n})}(x_0)\right|\leq \frac{1}{2}|\mathcal{Q}_{\tau_1 \theta_2 d_0/(2\sqrt{n})}(x_0)|.
\]
We conclude that
\begin{align*}
|A\cap \mathcal{Q}_{1}|\geq |A\cap \mathcal{Q}_{\tau_1 \theta_2 d_0/(2\sqrt{n})}(x_0)|\geq \frac{1}{2}|\mathcal{Q}_{\tau_1 \theta_2 d_0/(2\sqrt{n})}(x_0)|.
\end{align*}
This inequality contradicts \eqref{he1} because we took $t_0$ satisfying \eqref{he8}.
\end{proof}

\section{Examples}\label{sect8}
In this section, we provide an example that shows that Theorem \ref{mainthm2} fails if we assume only (A) but not \eqref{B}.

For $0<\beta_1<\beta_2 \leq2$ and $\beta_2<n$, it suffice to construct $\sigma_\epsilon(x):B_{2\sqrt{n}}\to[\beta_1,\beta_2]$ and $u_\epsilon$ for $\epsilon>0$, which is a non-negative solution of
\begin{align}
&\quad-\frac{\A(n,\sigma_\epsilon(x))}{\eta(n,\sigma_\epsilon(x))}(-\Delta)^{\sigma_\epsilon(x)/2}u_\epsilon(x)\nonumber\\
&:=\frac{\A(n,\sigma_\epsilon(x))}{\eta(n,\sigma_\epsilon(x))}\int_{\R^n}\delta(u_\epsilon,x,y)\frac{\eta(n,\sigma_\epsilon(x))dy}{|y|^{n+\sigma_\epsilon(x)}}=0\quad\mbox{in }B_{2\sqrt{n}},\label{eq82}
\end{align}
and we will see 
\begin{equation}\label{eq81}
\frac{\sup_{B_{1/4}}u_\epsilon}{\inf_{\mathcal{Q}_6}u_\epsilon}\to\infty\quad\mbox{as }\epsilon\to0.
\end{equation}
Here $\A(n,\sigma)>0$ is as in Example \ref{ex16} and $\eta(n,\sigma)>0$ is chosen so that the Fourier symbol of $(-\Delta)^{\sigma/2}$ is $-|\xi|^\sigma$ i.e. $\mathcal{F}[(-\Delta)^{\sigma/2}u](\xi)=-|\xi|^\sigma \mathcal{F}[u](\xi)$ for the Fourier transform $\mathcal{F}$.
Also note that we need to multiply $\A(n,\sigma_\epsilon(x))/\eta(n,\sigma_\epsilon(x))$ so that the nonlocal operator satisfies \eqref{A2}.
Let $\phi\in C^\infty_0(B_1)$ be a non-negative mollifier and $\phi_\epsilon(x):=\epsilon^{-n}\phi(x/\epsilon)$.
For $\beta_3\in(\beta_2,n)$, $l_0\in(0,1)$ and $N_1\in\N$, set
\[
u_\epsilon(x):=
\left\{
\begin{array}{cl}
(-\Delta)^{-\beta_3/2}\phi_\epsilon(x)&\mbox{for }x\in B_{3\sqrt{n}},\\
(-\Delta)^{-\beta_3/2}\phi_\epsilon(x)+l_0((|x|-3\sqrt{n})^+\wedge N_1)^{\beta_1}&\mbox{for }x\in \R^n\setminus B_{3\sqrt{n}},
\end{array}
\right.
\]
where $(-\Delta)^{-\beta_3/2}\phi_\epsilon(x)$ is the Riesz potential of $\phi_\epsilon$ defined by
\[
(-\Delta)^{-\beta_3/2}\phi_\epsilon(x):=c_{\beta_3}\int_{\R^n}\frac{\phi_\epsilon(x-y)}{|y|^{n+\beta_3}}dy
\]
for $c_{\beta_3}=\pi^{-n/2}2^{-\beta_3}\Gamma\left(\left(n-\beta_3\right)/2\right)/\Gamma\left(\beta_3/2\right) $.
For $x\in B_{2\sqrt{n}}$, we compute
\begin{align*}
&\quad-(-\Delta)^{\beta_2/2}u_\epsilon(x)\\
&=-(-\Delta)^{-(\beta_3-\beta_2)/2}\phi_\epsilon(x)+\int_{\R^n}l_0((|x+y|-3\sqrt{n})^+\wedge N_1)^{\beta_1}\frac{2\eta(n,\beta_2)dy}{|y|^{n+\beta_2}}\\
&\leq -(-\Delta)^{-(\beta_3-\beta_2)/2}\phi_\epsilon(x)+\int_{\R^n\setminus B_{\sqrt{n}}}l_0(|y|-\sqrt{n})^{\beta_1}\frac{2\eta(n,\beta_2)dy}{|y|^{n+\beta_2}}\\
&= -(-\Delta)^{-(\beta_3-\beta_2)/2}\phi_\epsilon(x)+C(n,\beta_1,\beta_2)l_0.
\end{align*}
As we have
\begin{align*}
(-\Delta)^{-(\beta_3-\beta_2)/2}\phi_\epsilon(x)
&=c_{(\beta_3-\beta_2)}\int_{\R^n}\frac{\phi_\epsilon(x-y)}{|y|^{n+\beta_3-\beta_2}}dy\\
&\geq c_{(\beta_3-\beta_2)}|2\sqrt{n}+1|^{-n+\beta_3-\beta_2}\int_{\R^n}\phi_\epsilon(x-y)dy\\
&= c_{(\beta_3-\beta_2)}|2\sqrt{n}+1|^{-n+\beta_3-\beta_2}>0
\end{align*}
for $x\in B_{2\sqrt{n}}$, there exists $l_0=l_0(n,\beta_1,\beta_2,\beta_3)>0$ such that
\begin{equation}\label{eq83}
-(-\Delta)^{-(\beta_3-\beta_2)/2}\phi_\epsilon(x)\leq0 \quad\mbox{for }x\in B_{2\sqrt{n}}.
\end{equation}
On the other hand, we have for $x\in B_{2\sqrt{n}}$,
\begin{align*}
&\quad-(-\Delta)^{\beta_1/2}u_\epsilon(x)\\
&=-(-\Delta)^{(\beta_3-\beta_1)/2}\phi_\epsilon(x)+\int_{\R^n}l_0((|x+y|-3\sqrt{n})^+\wedge N_1)^{\beta_1}\frac{2\eta(n,\beta_1)dy}{|y|^{n+\beta_1}}\\
&\geq -(-\Delta)^{(\beta_3-\beta_1)/2}\phi_\epsilon(x)+\int_{\R^n\setminus B_{5\sqrt{n}}}l_0((|y|-5\sqrt{n})^+\wedge N_1)^{\beta_1}\frac{2\eta(n,\beta_1)dy}{|y|^{n+\beta_1}}
\end{align*}
and
\begin{align*}
(-\Delta)^{(\beta_3-\beta_1)/2}\phi_\epsilon(x)
&=c_{(\beta_3-\beta_1)}\int_{\R^n}\frac{\phi_\epsilon(x-y)}{|y|^{n+\beta_3-\beta_1}}dy\\
&\leq c_{(\beta_3-\beta_1)}\|\phi_\epsilon\|_\infty\int_{B_{2\sqrt{n}+1}}\frac{dy}{|y|^{n+\beta_3-\beta_1}}dy
= C(n,\beta_1,\beta_3)\epsilon^{-n}.
\end{align*}
Since
\[
\lim_{N\to\infty}\int_{\R^n\setminus B_{5\sqrt{n}}}l_0((|y|-5\sqrt{n})^+\wedge N_1)^{\beta_1}\frac{2\eta(n,\beta_1)dy}{|y|^{n+\beta_1}}=\infty,
\]
there exists $N_1=N_1(n,\beta_1,\beta_2,\beta_3,\epsilon)$ such that
\begin{equation}\label{eq84}
-(-\Delta)^{\beta_1/2}u_\epsilon(x)\geq0\quad\quad\mbox{for }x\in B_{2\sqrt{n}}.
\end{equation}
From \eqref{eq83} and \eqref{eq84}, for each $x\in B_{2\sqrt{n}}$, there exists $\sigma_\epsilon(x)\in [\beta_1,\beta_2]$ such that \eqref{eq82} holds.
Moreover, because $u_\epsilon(x)=(-\Delta)^{-\beta_3/n}\phi_\epsilon(x)\to c_{\beta_3}|x|^{-n+\beta_3}$ for $x\in B_{2\sqrt{n}}$, we obtain \eqref{eq81}.


\subsection*{Acknowledgments.} The author wishes to express his thanks to Prof. Shigeaki Koike for many helpful suggestions during the preparation of the paper.
S. Kitano is supported by Grant-in-Aid for JSPS Fellows 21J10020.

\end{document}